\documentclass{article}

\usepackage{amsmath}
\usepackage{amssymb}
\usepackage{bm}
\usepackage{multicol}
\usepackage{multirow}
\usepackage{longtable}
\usepackage{pdflscape}
\usepackage{rotating}
\usepackage[sorting=none]{biblatex}
\usepackage{blindtext}
\usepackage{ragged2e}
\usepackage{graphicx}
\usepackage{titlesec}
\usepackage{subcaption}
\usepackage{csquotes}
\usepackage{float}
\usepackage[toc,title]{appendix}
\usepackage[linesnumbered,ruled,vlined]{algorithm2e}
\usepackage{hyperref}
\usepackage{cleveref}
\usepackage[acronym,style=super,nogroupskip,nonumberlist,automake]{glossaries-extra}

\makeglossaries
\setabbreviationstyle[acronym]{long-short}
\glsdisablehyper
\GlsXtrEnableEntryCounting{acronym}{1}

\DeclareMathOperator*{\argmin}{arg\,min}

\newacronym{ac}{AC}{alternating current}
\newacronym{acopf}{AC-OPF}{alternating current optimal power flow}
\newacronym{acpf}{AC-PF}{alternating current power flow}
\newacronym{admm}{ADMM}{alternating direction method of multipliers}
\newacronym{aladin}{ALADIN}{augmented Lagrangian alternating direction inexact Newton}
\newacronym[shortplural={BESSs},longplural={battery energy storage systems}]{bess}{BESS}{battery energy storage system}
\newacronym{bim}{BIM}{bus injection method}
\newacronym{capex}{CAPEX}{capital expenditure}
\newacronym{cop}{CoP}{coefficient of performance}
\newacronym{crf}{CRF}{capital recovery factor}
\newacronym{dcopf}{DC-OPF}{direct current optimal power flow}
\newacronym{dcpf}{DC-PF}{direct current power flow}
\newacronym[shortplural={DERs},longplural={distributed energy resources}]{der}{DER}{distributed energy resource}
\newacronym{des}{DES}{distributed energy system}
\newacronym{dod}{DoD}{depth of discharge}
\newacronym[shortplural={DSOs},longplural={distribution system operators}]{dso}{DSO}{distribution system operator}
\newacronym[shortplural={EVs},longplural={electric vehicles}]{ev}{EV}{electric vehicle}
\newacronym{elvtf}{ELV-TF}{European low voltage test feeder}
\newacronym{epsrc}{EPSRC}{Engineering and Physical Sciences Research Council}
\newacronym{hv}{HV}{high voltage}
\newacronym{ieee}{IEEE}{Institute of Electrical and Electronics Engineers}
\newacronym{lv}{LV}{low voltage}
\newacronym{milp}{MILP}{mixed-integer linear program}
\newacronym{minlp}{MINLP}{mixed-integer nonlinear program}
\newacronym{nlp}{NLP}{nonlinear program}
\newacronym{opex}{OPEX}{operational expenditure}
\newacronym{opf}{OPF}{optimal power flow}
\newacronym{pf}{PF}{power flow}
\newacronym[shortplural={PVs},longplural={photovoltaics}]{pv}{PV}{photovoltaic}
\newacronym[shortplural={SDGs},longplural={sustainable development goals}]{sdg}{SDG}{sustainable development goal}
\newacronym{sdp}{SDP}{semi-definite programming}
\newacronym{seg}{SEG}{smart export guarantee}
\newacronym{soc}{SoC}{state of charge}
\newacronym{socp}{SOCP}{second-order cone programming}
\newacronym{tac}{TAC}{total annualised cost}
\newacronym[shortplural={TSOs},longplural={transmission system operators}]{tso}{TSO}{transmission system operator}
\newacronym{un}{UN}{United Nations}
\newacronym{ved}{VED}{volumetric energy density}

\glssetcategoryattribute{acronym}{glossdesc}{title}
\MFUhyphentrue
\MFUnocap{on}
\MFUnocap{of}
\MFUnocap{the}
\MFUnocap{and}

\title{Distributed Energy System Design including Unbalanced AC Power Flow for Large LV Networks with ADMM}
\author{Robert Steven \and Oleksiy V. Klymenko \and Michael Short\\
  \texttt{m.short@surrey.ac.uk}}

\begin{document}
    \maketitle
    \begin{abstract}
With the addition of large numbers of \glspl{der} to distribution networks comes the increasing risk that their operation may violate the safety constraints of these networks. The problem considered in this paper is that of combined siting, sizing and dispatch of these \glspl{der}, also known as \gls{des} design, to help meet electrical and heat loads within the network. Here, the operation of these \glspl{der} is modelled, along with the unbalanced three-phase \gls{ac} power flow in the network. When this network power flow is considered, this admits a non-convex \gls{minlp} model formulation which scales poorly with network size in terms of solve time. To address this, the problem is decomposed into a series of  algorithmic steps, starting with a \gls{milp} formulation that does not consider network constraints, then fixing binary variables, adding power flow constraints and solving as a \gls{nlp} and finally removing operational binary variables and replacing them with a complementarity reformulation. As the main contributors to the overall solve time, the \gls{nlp} and Complementarity steps are solved using a hybrid spatial/temporal decomposition strategy and the \gls{admm} distributed optimisation method. Results are presented for networks based on the \glsxtrlong{elvtf} with up to $55$ loads and $120$ timepoints, with the \gls{admm} approach showing speed-ups of up to $13$x when considering parallel computation of the subproblems, for a maximum observed optimality gap of $0.61\%$.
    \end{abstract}
    
    \glsresetall

    \printglossary[type=\acronymtype]
    
    \section{Introduction}
The \Gls{un} \glspl{sdg} mandate moving to a net-zero carbon society, limiting the impacts of climate change and helping people around the world to gain fair and equitable access to energy \cite{UN_DESA}. Doing so requires widespread deployment of \glspl{der}, examples of which include \gls{pv} panels for power generation and \glspl{bess} for storage, as well as low-carbon heating solutions, for example heat pumps \cite{Pathak2022}, many of which will be sited within \gls{lv} distribution networks. These networks have traditionally been viewed as mainly passive, with the sole purpose of transporting electricity unidirectionally to users from the \gls{hv} transmission network, via a network containing electrical lines and one or more step-down transformers \cite{Kersting2012}. In distribution networks, the ratio of resistance to reactance (R/X) for a given line is typically high \cite{Dubey2023,Sereeter2017}, such that sinking/sourcing active power by end-users can cause significant fluctuations in voltage magnitude, which must be kept within certain limits to ensure the network is able to operate safely \cite{Frank2016}. \\
The large-scale inclusion of \glspl{der} and low-carbon heating in these networks presents a number of challenges in this respect. The ability of end-users to export power back into the network via schemes such as the \gls{seg} in the UK \cite{SEG_Ofgem} leads to bidirectional power flow, with associated changes to bus voltages within the network. Low-carbon heating technologies, for example heat pumps, may also place significant additional electrical demand on the network if building heating loads were previously being met by natural gas-based technologies such as condensing boilers. Therefore, if these \gls{der} and low-carbon heating technologies are installed and operated without consideration of the underlying network they are placed within, which may not have been initially constructed with their use in mind, then the possibility of network safety violations occurring increases \cite{Steven2025b}. It may therefore be a benefit for network operators to determine the effect that the cost-optimal installation of these technologies, in terms of capacity and location, as well as their active control through the use of \gls{opf}, will have on their network, to mitigate the propensity for the violation of operating conditions. \\
\Gls{des} models that include the sizing, siting and dispatch of \glspl{der} within a network can be used to address this, with a review by De Mel et al. \cite{DeMel2022a} detailing a range of different modelling formulations. Where these models take into account the constraints of the underlying network, they incorporate \gls{pf} constraints used in \gls{opf} models (see Frank \& Rebennack \cite{Frank2016} for a detailed introduction to \gls{opf} models) and must remain feasible with respect to network operating constraints, such as voltage magnitude upper and lower bounds \cite{DeMel2024CompReform}. Where the full \gls{acpf} constraints are considered (\gls{acopf}), these add non-convex nonlinearity to the model \cite{Frank2016}. Furthermore, unlike the \gls{hv} transmission network, the power flows in each phase of a distribution network cannot be assumed to be balanced \cite{Dubey2023}. This requires the explicit consideration of each phase (and the corresponding inter-phase coupling), which adds additional modelling complexity \cite{DeMel2024CompReform}. \\
Solving \gls{des} models entails finding the optimal set of installations of \glspl{der} and their corresponding dispatches, with \enquote{optimal} typically being measured using an environmental or economic objective function \cite{DeMel2022a}. Approaches based on mathematical optimisation can give provably optimal results and have been used extensively for this purpose \cite{DeMel2022a}, but can struggle as a result of power flow constraints present in the model adding non-convex nonlinearity, especially for larger problem sizes and when mixed-integer variables are also included \cite{DeMel2024CompReform}. Additional approaches that address this include metaheuristic algorithms and convex relaxations. Numerous metaheuristic algorithms are present in the literature, however in general they contain a number of drawbacks, including that they cannot guarantee convergence to an optimal solution \cite{DeMel2022a}, and they may suffer from long compute times \cite{Papadimitrakis2021}. Convex relaxations such as \gls{sdp} \cite{Lavaei2012,Gan2014} and \gls{socp} \cite{Bobo2021} have been successfully applied to \gls{opf} formulations and have demonstrated exact relaxation under certain conditions \cite{Dubey2023}. These may not apply to unbalanced three-phase formulations however \cite{Dubey2023} and whilst their convexity may improve the tractability of the problem, relaxations such as \gls{sdp} introduce additional variables which can become significant for larger problem sizes \cite{Erseghe2015}. \\
When considering large-scale optimisation problems in power systems, the use of distributed optimisation can ensure that the computation of the solution remains tractable by taking the main problem and decomposing it into several smaller subproblems \cite{Molzahn2017}. These subproblems are then solved separately, whilst the distributed optimisation algorithm must also ensure that their solutions remain feasible with respect to the main problem \cite{Boyd2011}. The \gls{opf} problem has been tackled by a number of distributed optimisation approaches \cite{Molzahn2017,Kargarian2018,Alkhraijah2024,Steven2025c}, from which the \gls{admm} \cite{Boyd2011} has proven to be particularly popular \cite{Molzahn2017}. In general, \gls{admm} requires convex formulations for proof of convergence \cite{Boyd2011}, such as \gls{sdp}/\gls{socp} formulations or linearisations such as \gls{dcopf} \cite{Frank2016} or LinDistFlow \cite{Gan2014}. Further modifications to this approach address this directly, including an \gls{admm}-based technique proposed by Erseghe \cite{Erseghe2015} and the \glsxtrfull{aladin} approach by Engelmann et al. \cite{Engelmann2019}. The approach by Erseghe \cite{Erseghe2015} is tailored for use with non-convex \gls{acpf} constraints, with a proof of convergence that assumes the solver in each subproblem can return a locally optimal solution. This has been applied to large-scale \gls{acopf} problems, demonstrating a speed-up in determining the solution compared to a centralised problem formulation, by Guo et al. \cite{Guo2017}. \\
As well as \gls{opf} problems, a number of approaches have applied \gls{admm} specifically to the problem of sizing \& siting \glspl{der} within a network, in a variety of different ways. Nick et al. \cite{Nick2015} present a formulation for calculating the optimal siting and sizing of \glspl{bess} to provide ancillary services to the distribution network they are installed in. The authors decompose the optimisation between investment and a number of operational scenarios, optimising the scenarios in parallel for a given investment decision. Their formulation includes a \gls{socp} convex relaxation of the \gls{acpf} constraints and they report significant speed-up in overall compute time from the use of \gls{admm}. A scenario-based decomposition is also presented by Liu et al. \cite{Liu2024}, using linearised \gls{pf} constraints and optimising battery storage over a number of scenarios and \gls{pv} under a worst-case scenario. Battery storage is a further focus of Sayfutdinov et al. \cite{Sayfutdinov2020}, who provide a formulation that includes a degradation model for the installed \gls{bess} in their problem that considers both siting \& sizing as well as technology selection. They use a mixed-integer convex formulation, including the \gls{dcpf} approximation for network power flow and decompose their main problem into subproblems for each bus and battery technology option. A previous work by the authors of this paper also examined \gls{pv} and \gls{bess} sizing and dispatch \cite{Steven2024}, using the \gls{dcpf} approximation and a temporal decomposition. The approach used by Liu et al. \cite{Liu2023} employs the full non-convex \gls{acpf} constraints in a distribution network planning problem, whereby the network is constructed and renewable generation installed within it by a number of participants, solved using \gls{admm} for the \gls{ieee} 33-bus distribution system. The \gls{pf} formulation they employ does not consider unbalanced three-phase operation however, which may limit its applicability. \\
This unbalanced operation is considered in the planning methods employed by both Shi et al. \cite{Shi2025} and Zhu et al. \cite{Zhu2026}, focussing on the interconnection with transport networks and the interaction of multiple distribution networks respectively. The method presented by Shi et al. \cite{Shi2025} aims to optimise the installation of \gls{ev} chargers, alongside upgrades to the distribution and transport networks themselves, where they note the propensity for the additional loads presented by \glspl{ev} to worsen voltage imbalance within the network. This is solved using a \gls{milp} formulation however, with the authors presenting a linearisation procedure for their representation of both power and transport networks. These networks form their decomposed subproblems, which they solve with an iterative process allowing them to interface \gls{admm} with binary variables. The decomposition used by Zhu et al. \cite{Zhu2026} is instead carried out at a regional level, where they present a formulation for multiple distribution networks connected together via soft open points coupled with \gls{bess}. They employ an \gls{sdp} convex relaxation, which they empirically demonstrate to be numerically exact. \\
\Gls{admm} has also been applied to adjacent problems, such as by Gu et al. \cite{Gu2024} who employ it with spatial decomposition, for use in a \glsxtrshort{tso}-\glsxtrshort{dso} (\glsxtrlong{tso} - \glsxtrlong{dso}) coordination problem to determine optimal \gls{der} installation by \glspl{dso} whilst taking into account their ability to support the upstream \glsxtrshort{tso}. The algorithm is also applied by Kilkki et al. \cite{Kilkki2018} to a formulation allowing for residential participation in energy markets, through the optimisation of energy storage and residential heating. The participation is aggregated, with the commitments of each individual resident forming the subproblems that \gls{admm} is decomposed on. \\
This existing works applying \gls{admm} to the sizing \& siting problem highlights a number of approaches to decomposing the problem space, along with a variety of problem formulations and particular focusses, such as on \gls{bess} placement. This paper proposes to address a number of current omissions, by presenting an approach which does not rely on convex relaxations or linearisations, uses a full three-phase unbalanced power flow formulation and includes in its consideration the installation of heat pumps as a low-carbon heating option. The presented model combines optimal siting/sizing and dispatch for \glspl{der} and heating technologies within a distribution network, taking into account the unbalanced three-phase \gls{acpf} constraints within the distribution network. The problem decomposition is a hybrid spatial/temporal approach, solved using the \gls{admm}-based method presented by Erseghe \cite{Erseghe2015}. This paper adapts the given approach to the multi-timepoint formulation required for the combined siting/sizing and dispatch model, with details given for the decomposition structure and partitioning methods employed. It is envisaged that the results can help examine the feasibility of installing \glspl{der} and low-carbon heating technologies for a given network, both in terms of their location/sizing as well as their operation. The computational method presented further alleviates the burden of solving these problems in larger networks, allowing them to be computed in reasonable time under the assumption of parallel computation of the decomposed subproblems.
    
    \section{Methodology}
\subsection{Nomenclature}
    The problem considered here pertains to the installation of \glspl{der} within an electrical distribution network, denoted as an undirected graph $\mathcal{G}=\{\mathcal{N},\mathcal{E}\}$ containing nodes and edges, referred to respectively as buses ($n=\{1,2,\ldots,N\}\in\mathcal{N}$) and branches connecting bus $i$ to $j$ ($(i,j)\in\mathcal{E}$). A subset of these buses contains loads, $\mathcal{N}_L\subseteq\mathcal{N}$, defined by the network formulation, for which there are both electrical and heating loads, $E^{load}$ and $H^{load}$. These $\mathcal{N}_L$ buses are considered as candidates for the installation of \glspl{der} and heating technologies to assist in meeting these loads. 
    A single bus within the network is defined as the \enquote{slack} bus \cite{Frank2016}, modelling the connection to the wider power network, which has a fixed reference voltage magnitude and phase, and allows unconstrained real and reactive power sink/sourcing. Buses and branches within the network can each have up to $3$ phases, defined as $\psi\in\phi=\{A,B,C\}$. Each bus has a corresponding complex voltage phasor for each phase, expressed in polar coordinates with magnitude $V^{\psi}_i$ and angle $\theta^{\psi}_i$ for phase $\psi$ and bus index $i$. Buses containing loads, and the slack bus, may also have real and reactive power injections into the network, defined as $P^{\psi}_i$ and $Q^{\psi}_i$ respectively \cite{Frank2016}, with positive values indicating injection into the network. Branches define connections between buses within the network, with each connection modelled as a $3\times 3$ (for three phases) complex admittance matrix, $Y_{i,j}$ \cite{Kersting2012}, for the connection between buses $i$ and $j$, where off-diagonal terms signify electrical coupling between phases. Each complex admittance matrix, $Y_{i,j}$, is made up of the conductance (real), $G_{i,j}$, and susceptance (imaginary) matrices $B_{i,j}$.

\subsection{DES Design Model}
    The \gls{des} design model considered here is built on the one presented by De Mel et al. \cite{DeMel2024CompReform} with the technologies that can be installed at each load presented in Table \ref{tab:der_technologies}.

    \begin{table}[H]
        \centering
        \caption{\Glsxtrshort{der} and Low-Carbon Heating Technologies}
        \begin{tabular}{|l|l|}
            \hline
            \textbf{Technology} & \textbf{Output} \\
            \hline
            \Gls{pv} & Electrical energy generated \\
            \hline
            Battery storage & Electrical energy charged, discharged and stored \\
            \hline
            Boilers & Heat energy generated \\
            \hline
            Heat Pumps & Heat energy generated \\
            \hline
            Hot Water Tanks & Heat energy charged, discharged and stored \\
            \hline
        \end{tabular}
        \label{tab:der_technologies}
    \end{table}

    Modelling equations for these technologies are detailed by De Mel et al. \cite{DeMel2024CompReform,DeMel2023HeatDecarb} and presented in full in Appendix \ref{app:modelling_components}. Constraints within the model cover the operation of each technology, their integration to assist in meeting the electrical and heating loads and linking constraints enforcing that the same installation decisions are made across all timepoints. \\
    The model has an economic objective function, referred to as \gls{tac}, made up of \gls{capex} and \gls{opex} components pertaining to installations and operations within the network. The \gls{crf} is then used to consider the payback over the total number of years considered, taken here to be $20$. \\
    The model itself considers a $24$h period for each of the four seasons, with a resolution of $1$h, by averaging values across each season. Input parameters such as building electrical/heat loads, ambient temperature and irradiance are then provided for each of these timepoints. Further parameters covering the \gls{der} technologies and scalar model parameters can be found in Appendix \ref{app:parameters_variables}.

    \subsubsection{Model Formulation}
        Determining the optimal combined siting/sizing and dispatch of \glspl{der} admits a \gls{minlp} formulation, with the following model components defined in Table \ref{tab:model_class_table}. From the table, it can be seen that the installation and operation of the \glspl{der} has a mixed-integer linear form, whilst constraints governing power flow within the network add non-convex nonlinearity to the model.

		\begin{table}[H]
            \centering
            \caption{Model Components}
            \begin{tabular}{|l|l|p{45mm}|}
                \hline
                \textbf{Component} & \textbf{Class} & \textbf{Description} \\
                \hline
                Objective Function & Linear & Combined \gls{capex} and \gls{opex} expenditure \\
                \hline
                Installation Decisions (Siting) & Integer/Binary & Decision to install \gls{der} at a given bus \\
                \hline
                \Gls{der} Sizing \& Operation & Linear & Constraints governing \gls{der} operation \\
                \hline
                Dispatch & Binary & Battery charge/discharge and grid buy/sell \\
                \hline
                Grid Power Flow & Nonlinear & Non-convex \gls{acpf} constraints \\ 
                \hline
            \end{tabular}
            \label{tab:model_class_table}
		\end{table}

    \subsubsection{Solution Strategy}
\begin{figure}[H]
            \centering
            \includegraphics[width=0.85\linewidth]{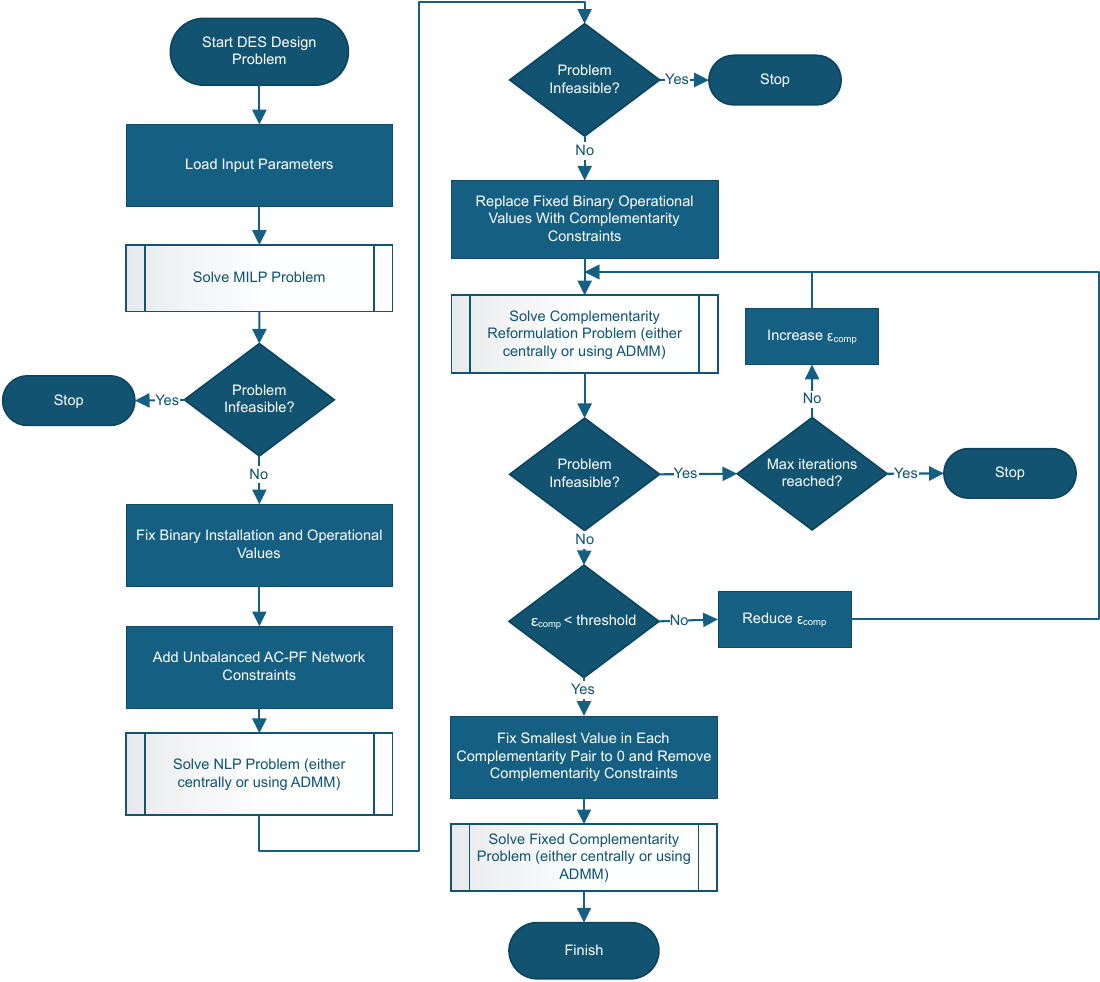}
            \caption{Solution Algorithm Overview}
            \label{fig:solution_algorithm_overview}
        \end{figure}

        The overall solution strategy, similar to the approach employed by De Mel et al. \cite{DeMel2024CompReform}, is shown in Algorithm \ref{alg:central_solve} and further illustrated in Figure \ref{fig:solution_algorithm_overview}. The problem is solved in three main stages, as an \gls{milp} without \gls{pf} constraints, as an \gls{nlp} where binary variables are fixed and nonlinear \gls{acpf} constraints added and finally with complementarity reformulations replacing the operational binary variables. Here, the \gls{nlp} and complementarity reformulations are either solved centrally or decomposed and solved using distributed optimisation (\gls{admm}).

        \begin{algorithm}[H]\label{alg:central_solve}
            Solve as \glsxtrshort{milp}\label{line:central_milp_solve}\;
            Binary/integer decision values are fixed and \glsxtrshort{acpf} constraints added\;
            Model is re-solved as an \glsxtrshort{nlp} to determine \glsxtrshort{der} sizing values that respect \glsxtrshort{acpf} constraints\label{line:central_nlp_solve}\;
            Operational binary decisions are un-fixed and replaced with complementarity reformulations to allow these to influence the optimisation\label{line:central_comp_solve}\;
            \caption{Model Central Solve Algorithm}
        \end{algorithm}

    \subsubsection{Robust Season}
        In addition to the representative $24$h periods for each season that are considered, a \enquote{robust} $24$h period is also included, giving a total of $120$ timepoints. To include this in the main model, linking constraints for the siting/sizing decision variables are added, such that these decisions must meet the robust loads and environmental conditions. The \gls{opex} cost for the robust season is not considered however. In this way, the model seeks to make \gls{capex} decisions that are robust to the input load and environmental parameters whilst minimising for the average operational cost of the technologies. Further details for how the robust season values are generated can be found in Section \ref{building_load_environ}.
	
    \subsubsection{Heat Pump Model}
        The main equations for the modelling of heat pumps are given in De Mel et al. \cite{DeMel2023HeatDecarb}, however the method for determining the \gls{cop} and the capacity of each heat pump has been altered. Instead of a piecewise linearisation, these are now derived by fitting a sigmoid function to the data-points provided by the manufacturer \cite{Mitsubishi_i} to improve the extrapolation of the parameters to ambient temperatures outside of the range provided. The following logistic function is used:
		\begin{equation}
            \frac{Le^{k\left(x-x_0\right)}}{1+e^{k\left(x-x_0\right)}}+c
		\end{equation}
		With the ambient temperature ($T^{amb}$) as the input, $x$, and the curve fitting to manufacturer-provided values carried out using the nonlinear least squares curve\_fit() function in SciPy \cite{SciPy,Vugrin2007}. The fit for this was seen to be satisfactory for the provided heat pumps.
 
    \subsubsection{Power Flow Model}\label{power_flow_model}
		The \gls{des} design model, not considering power flow within the network, is classified as a \gls{milp} model (as can be seen in Table \ref{tab:model_class_table}), which can be efficiently solved at scale with existing solvers such as Gurobi \cite{Gurobi12}. To capture the flow of power within the network however, presents well-known additional complexity. To model the flow of \gls{ac} power within the network requires the use of non-convex, nonlinear \gls{acpf} constraints \cite{Frank2016}. Solving this in direct combination with the \gls{milp} \gls{des} design model results in an \gls{minlp} formulation which is known to scale poorly in terms of compute time for the problem size (the number of variables and constraints contained in the model). This was addressed by De Mel et al. \cite{DeMel2024CompReform} using a three-stage formulation whereby the model is solved as an \gls{milp} with the \gls{dcpf} approximation of the power flow constraints to preserve linearity, then the binary/integer decisions are fixed and the \gls{dcpf} constraints are replaced with the full three-phase \gls{acpf} \gls{bim} constraints and solved as a non-convex \gls{nlp}.

        \paragraph{MILP Model}
            A similar approach to the one outlined in Section \ref{power_flow_model} is used here, except the \gls{milp} problem is solved without \gls{dcpf} constraints, resulting in a formulation which does not initially capture interactions with the underlying network. In this way, buses which are installation candidates are able to install technologies and accordingly import/export power from/to the grid freely. This then provides a lower bound for the \gls{tac} of the model which may be infeasible with respect to the \gls{acpf} constraints, from which subsequent algorithm steps then work to improve. The use of a linearised \gls{pf} model, namely LinDistFlow (presented for three phase networks by Gan \& Low \cite{Gan2014}), was considered here but ultimately not included as the likely higher initial lower bound may have prematurely discarded subsequent feasible results in later algorithm stages.

		\paragraph{AC BIM PF}
            For the \gls{nlp} problem, the \gls{milp} solution (including \gls{der} sizing, siting \& dispatch and bus power injections (real, $P^{\phi}_j$ and reactive, $Q^{\phi}_j$)) is used as warm-start values and \gls{acpf} constraints are added, using the following \gls{bim} equations \cite{DeMel2024CompReform}:
            \begin{align}
                P^{\phi}_i=V^{\phi}_i \sum_{i} \left[ \sum_{\psi} \left[ V^{\psi}_j \left( G^{\phi\psi}_{i,j}cos\left( \theta^{\phi}_i-\theta^{\psi}_j \right) + B^{\phi\psi}_{i,j}sin\left( \theta^{\phi}_i-\theta^{\psi}_j \right) \right) \right] \right], \forall i\in\mathcal{N} \label{eq:bim_real_power_flow} \\
                Q^{\phi,}_j=V^{\phi}_i \sum_{i} \left[ \sum_{\psi} \left[ V^{\psi}_j \left( G^{\phi\psi}_{i,j}sin\left( \theta^{\phi}_i-\theta^{\psi}_j \right) - B^{\phi\psi}_{i,j}cos\left( \theta^{\phi}_i-\theta^{\psi}_j \right) \right) \right] \right], \forall i\in\mathcal{N} \label{eq:bim_react_power_flow}
            \end{align}
            \begin{align*}
                \text{where:} \\
                V^{\phi}_i &= \text{\RaggedRight Bus voltage magnitude for bus i and phase $\phi$} \\
                \theta^{\phi}_i &= \text{\RaggedRight Bus voltage angle for bus i and phase $\phi$} \\
                G^{\phi\psi}_{i,j} &= \text{\RaggedRight Branch conductance for branch (i,j) and pair of phases ($\phi$, $\psi$}) \\
                B^{\phi\psi}_{i,j} &= \text{\RaggedRight Branch susceptance for branch (i,j) and pair of phases ($\phi$, $\psi$}) \\
            \end{align*}
            with:
            \begin{align}
                Y^{\phi\psi}_{i,j}=G^{\phi\psi}_{i,j}+\bm{j}B^{\phi\psi}_{i,j}
            \end{align}
            \begin{align*}
                \text{where:} \\
                Y^{\phi\psi}_{i,j} &= \text{\RaggedRight Nodal admittance for branch $(i,j)$ and pair of phases ($\phi$ and $\psi$}) \\
            \end{align*}
            In the same way as De Mel et al. \cite{DeMel2024CompReform}, branches were modelled using the Approximate Line Segment Model \cite{Kersting2012} and transformers using the models presented by Sereeter et al. \cite{Sereeter2017}, based on those presented by Chen et al. \cite{Chen1991}. External grid source modelling uses the approach given by OpenDSS \cite{opendss}.

        \paragraph{AC Branch PF}
            Where branch power flows are also considered, these are defined as \cite{Gan2014}:
            \begin{align}
                \bm{Pb}_{\left(i,j\right)}=Re\left(diag\left(\tilde{\bm{V}}_{i}\left(Y_{\left(i,j\right)}\left(\tilde{\bm{V}}_{i}-\tilde{\bm{V}}_{j}\right)\right)^{H}\right)\right), \forall (i,j)\in\mathcal{E} \label{eq:bfm_real_branch_power_flow} \\
                \bm{Qb}_{\left(i,j\right)}=Im\left(diag\left(\tilde{\bm{V}}_{i}\left(Y_{\left(i,j\right)}\left(\tilde{\bm{V}}_{i}-\tilde{\bm{V}}_{j}\right)\right)^{H}\right)\right), \forall (i,j)\in\mathcal{E} \label{eq:bfm_imag_branch_power_flow}
            \end{align}
            \begin{align*}
                \text{where:} \\
                \bm{Pb}_{\left(i,j\right)} &= \text{\RaggedRight Real power flow for branch (i,j)} \\
                \bm{Qb}_{\left(i,j\right)} &= \text{\RaggedRight Imaginary power flow for branch (i,j)} \\
                \tilde{\bm{V}}_i &= \text{\RaggedRight Vector of voltage phasors, $\tilde{V}=V\angle\theta$, for each phase, for bus $i$} \\
                {\{\cdot\}}^H &= \text{\RaggedRight Hermitian transpose}
            \end{align*}

    \subsubsection{Complementarity Reformulations}
        To avoid admitting an overly conservative solution due to the inflexibility of the fixed operational binary decisions on battery and grid dispatch, the approach by De Mel et al. \cite{DeMel2024CompReform} added complementarity reformulations to these operational decisions (purchase/sale of grid power and battery charge/discharge dispatch decisions). Due to the increased nonlinearity introduced by the complementarity constraints, the \gls{nlp} solution (step \ref{line:central_nlp_solve} in Algorithm \ref{alg:central_solve}) is used as a good warm-start point. In this way, a trade-off is made whereby the compute time for larger models is reduced, rendering the model tractable, whilst still incorporating binary installation decisions and allowing binary dispatch decisions to affect the model outcome in conjunction with the \gls{acpf} constraints. \\
        The constraints that are replaced are as follows:
        \begin{align}
            E^{batt,charge}_{h,t,c} \leq M^{batt,chg} \cdot Q^{batt}_{h,t,c}, \forall h\in\mathcal{H}, t\in\mathcal{T}, c\in\mathcal{C} \\
            E^{batt,discharge}_{h,t,c} \leq M^{batt,chg} \cdot \left(1 - Q^{batt}_{h,t,c} \right), \forall h\in\mathcal{H}, t\in\mathcal{T}, c\in\mathcal{C} \\
            E^{grid}_{h,t} \leq M^{grid} \cdot \left(1 - X_{h,t}\right), \forall h\in\mathcal{H}, t\in\mathcal{T} \\
            E^{PV,sold}_{h,t} \leq M^{grid} \cdot X_{h,t}, \forall h\in\mathcal{H}, t\in\mathcal{T}
        \end{align}
        \begin{align*}
            \text{where:} \\
            \mathcal{H} &= \text{\RaggedRight Set of all loads} \\
            \mathcal{T} &= \text{\RaggedRight Set of all timepoints} \\
            \mathcal{C} &= \text{\RaggedRight Set of all batteries}
        \end{align*}
        are replaced with:
        \begin{align}
            E^{batt,charge}_{h,t,c} \cdot E^{batt,discharge}_{h,t,c} \leq \epsilon^{complementarity}, \forall h\in\mathcal{H}, t\in\mathcal{T}, c\in\mathcal{C} \\
            E^{grid}_{h,t} \cdot E^{PV,sold}_{h,t} \leq \epsilon^{complementarity}, \forall h\in\mathcal{H}, t\in\mathcal{T}
        \end{align}
        where the algorithm presented by De Mel et al. \cite{DeMel2024CompReform} uses an iterative process that repeatedly decreases $\epsilon^{complementarity}$ until below a fixed threshold, at which point overall convergence is declared. Following this convergence, in each pair of complementarity variables the smaller value is fixed to $0$ and the problem is re-solved to recover the \enquote{true} complementarity for each variable pair.

\subsection{Solving Model Centrally}
    When solving centrally, the model is solved with the steps shown in Algorithm \ref{alg:central_solve}. By solving over a range of network sizes, the effect on solution optimality and solve time can be demonstrated. Long solve times for the \gls{nlp} and complementarity solves motivates the reformulation of the problem to be solved in a distributed manner, as detailed in Section \ref{admm_formulation}. Note that for the purposes of this study, the \gls{seg} tariff value of $0.04\ \pounds/kWh$ from \cite{OVO_SEG} was insufficient to show changes in the solution between steps \ref{line:central_milp_solve}, \ref{line:central_nlp_solve} and \ref{line:central_comp_solve} in Algorithm \ref{alg:central_solve}. Therefore, to demonstrate the utility of this approach, the value was raised to $0.132\ \pounds/kWh$.

\subsection{ADMM Problem Formulation}\label{admm_formulation}
    \Gls{admm} is an algorithm that allows for the distributed optimisation of a model, through a \enquote{decomposition-coordination} approach \cite[p. 4]{Boyd2011}. The approach can be illustrated (with equations from Boyd et al. \cite{Boyd2011}) as follows for a global consensus problem, shown in \cref{eq:admm_consensus_problem}.
    \begin{alignat}{2}\label{eq:admm_consensus_problem}
        min\quad & f(x)=\sum_{i=1}^N f_i\left(x\right)
    \end{alignat}
    This is transformed into a form with a separable objective function in \cref{eq:admm_consensus_problem_seperable}.
    \begin{alignat}{2}\label{eq:admm_consensus_problem_seperable}
        min\quad & \sum_{i=1}^N f_i\left(x_i\right) \\
        \text{s.t.}\quad & x_i-z=0,i=1,\ldots,N
    \end{alignat}
    Where each objective function component, $f_i$, is now operating on its own copy of the problem variables, with the added constraint that there must be consensus across all components on the value of this variable. \\
    This then gives the following augmented Lagrangian:
    \begin{equation}
         L_{\rho}\left(x_1,\ldots,x_n,z,\lambda\right) = \sum_{i=1}^N\left[f_i(x_i)+\lambda^T_i(x_i-z)+\left(\frac{\rho}{2}\right)\left\Vert x_i-z\right\Vert^2_2\right]\label{eq:admm_aug_lagrange}
    \end{equation}
    with the iterations of the algorithm as:
    \begin{align}
		x_i^{k+1}:=\argmin_{x_i}L_{\rho}\left(x_i,z^{k},\lambda^{k}_i\right)=\argmin_{x_i}\left(f_i(x_i)+\left(\lambda^k_i\right)^T(x_i-z^k)+\left(\frac{\rho}{2}\right)\left\Vert x_i-z^k\right\Vert^2_2\right),i=1,\ldots,N \label{eq:admm_x_upd} \\
		z^{k+1}:=\argmin_{z}L_{\rho}\left(x^{k+1}_i,z,\lambda^{k}_i\right)=\frac{1}{N}\sum_{i=1}^N\left[x^{k+1}_i+\left(\frac{1}{\rho}\right)\lambda_i^k\right] \label{eq:admm_z_upd} \\
		\lambda^{k+1}_i:=\lambda^{k}_i+\rho\left(x^{k+1}_i-z^{k+1}\right),i=1,\ldots,N \label{eq:admm_lambda_upd}
    \end{align}
    Each of the subproblems in \cref{eq:admm_x_upd} now only depends on its own copy of the problem variables, $x_i$, as well as $z^k$ and $\lambda^k_i$ which are held as fixed values until the following algorithm steps. Therefore, these subproblems can now be updated in a parallel, synchronous manner (\textbf{decomposition}) \cite{Boyd2011}, meaning that the algorithm cannot advance to \cref{eq:admm_z_upd} until all subproblems have completed. Following that, it can be seen that \cref{eq:admm_z_upd} and \cref{eq:admm_lambda_upd} can then be trivially updated (\textbf{coordination}). These iterations are repeated until a stopping criterion, such as the reduction of the primal residual to below a set threshold ($\left\Vert x^{k+1}_i-z^{k+1}\right\Vert_{\infty}\leq\epsilon$), is met. \\
    The \gls{admm}-based form used in this paper is derived from Erseghe \cite{Erseghe2015}, which is adapted from the main \gls{admm} method to give convergence guarantees for the non-convex \gls{acpf} constraints (the standard \gls{admm} method requires convexity for its proof of convergence \cite{Boyd2011}). This formulation decomposes the network spatially, creating a number of subproblems that also encompass the electrically connected neighbours of buses in the subproblem. The complicating constraints are then those governing power flow through the tie-line branches between buses in the subproblem and their electrically connected neighbours. To solve, consensus is driven on the voltage magnitudes and angles of these overlapping buses (which as can be seen in \cref{eq:bim_real_power_flow} and \cref{eq:bim_react_power_flow} (for unbalanced three-phase) govern power flow through a given branch) by exchanging information between neighbouring subproblems. Here, this approach is adapted for the multi-timepoint \gls{des} design model, with consensus additionally being driven on branch real and reactive power flows for the tie-line branches. The steps for solving with \gls{admm} are given in Algorithm \ref{alg:admm_solve} and the reader is referred to Erseghe \cite{Erseghe2015} for the full formulation applied to the \gls{acopf} problem.
	
    \subsubsection{ADMM Consensus}
		Consensus is driven on tie-line branches between partitioned subproblems. As such, the variables requiring consensus are the voltage magnitude ($\bm{V}_t$) and angle ($\bm{\theta}_t$) at each end of the branch (in both subproblems) and the branch real ($\bm{Pb}_t$) and reactive ($\bm{Qb}_t$) power flow. For load subproblems, state variables are therefore defined for each timepoint as:
		\begin{equation}
            \bm{x}_{t}=\begin{bmatrix}
                \bm{V}_t, & \bm{\theta}_t, & \bm{Pb}_t, & \bm{Qb}_t
            \end{bmatrix}^T
		\end{equation}
		with each variable then consisting of a column vector for all buses/branches and phases. For voltage magnitude for example:
		\begin{equation}
            \bm{V}_t = \begin{bmatrix}
                V^A_{1,t}, & \cdots & V^A_{n,t}, & V^B_{1,t}, & \cdots & V^B_{n,t}, & V^C_{1,t}, & \cdots & V^C_{n,t}
            \end{bmatrix}^T
		\end{equation}
		is a column vector containing all the voltage magnitudes for all buses in the subproblem, across all phases. Although this is a very large vector, it is then multiplied by a sparse $\bm{A}$ consensus matrix when used in \gls{admm}. \\
		The overall state variable for a given subproblem is then the vertical concatenation of all timepoints (shown for $n$ timepoints):
		\begin{equation}
            \bm{x}=\begin{bmatrix}
                \bm{x}_1, & \bm{x}_2, & \cdots, & \bm{x}_n
            \end{bmatrix}^T
		\end{equation}
		The state variables for no-load partitions are considered in the same way, but are not concatenated across timepoints as these are considered separately (see Section \ref{problem_decomposition}). In the no-load partitions, the power flow on phases of tie-line branches where no load is connected to that phase on the other side are fixed to zero. \\
		The consensus matrix, $\bm{A}_{i}$, for each subproblem is then constructed, with a single row per consensus constraint in each problem. By setting the corresponding column values to $1$, constraints of the following form for each tie-line branch, phase and timepoint can be constructed:
		\begin{align}
            V^{\phi}_{n,t} = V^{\phi}_{i,t}\label{eq:v1_consensus} \\
            V^{\phi}_{m,t} = V^{\phi}_{j,t}\label{eq:v2_consensus} \\
            \theta^{\phi}_{n,t} = \theta^{\phi}_{i,t}\label{eq:theta1_consensus} \\
            \theta^{\phi}_{m,t} = \theta^{\phi}_{i,t}\label{eq:theta2_consensus} \\
            Pb^{\phi}_{l,t} = Pb^{\phi}_{h,t} \\
            Qb^{\phi}_{l,t} = Qb^{\phi}_{h,t}
		\end{align}
		where the pairs $\left(n,i\right)$ denote the same bus at one end of the tie-line branch in the different subproblems, $\left(m,j\right)$ the bus at the other end and $\left(l,h\right)$ the branch itself. \\
        As presented by Erseghe \cite{Erseghe2015}, voltage-based constraints can be reformulated as an addition and subtraction:
        \begin{align}
            V^{\phi}_{n,t}+V^{\phi}_{m,t}=V^{\phi}_{i,t}+V^{\phi}_{j,t} \\
            V^{\phi}_{n,t}-V^{\phi}_{m,t}=V^{\phi}_{i,t}-V^{\phi}_{j,t}
        \end{align}
        which was also used here, in place of the straight voltage equality consensus constraints (\crefrange{eq:v1_consensus}{eq:theta2_consensus}). \\
        To improve the numerical performance of solving the parallel subproblems (\cref{eq:admm_x_upd}), a scaling factor was applied to the \gls{milp} \gls{des} cost components in the objective function. This was done by normalising the \gls{tac} of each subproblem to the value given by the the central \gls{milp} solve (line \ref{line:central_milp_solve} in Algorithm \ref{alg:central_solve}). For comparison purposes, this same approach was also applied to the centralised formulation.

        \begin{algorithm}[H]\label{alg:admm_solve}
            Solve centrally as \glsxtrshort{milp}\label{line:dist_milp_solve}\;
            Binary/integer decision values are fixed and \glsxtrshort{acpf} constraints added\;
            Model is decomposed into a number of subproblems\;
            Subproblems are solved as an \glsxtrshort{nlp} using \glsxtrshort{admm}\label{line:dist_nlp_solve}\;
            Operational binary decisions are un-fixed and replaced with complementarity reformulations and re-solved using \glsxtrshort{admm}\label{line:dist_comp_solve}\;
            \caption{Model Distributed Solve Algorithm}
        \end{algorithm}
	
    \subsubsection{Problem Decomposition}\label{problem_decomposition}
\begin{figure}[h]
            \centering
            \begin{subfigure}{0.35\textwidth}
                \includegraphics[width=\textwidth]{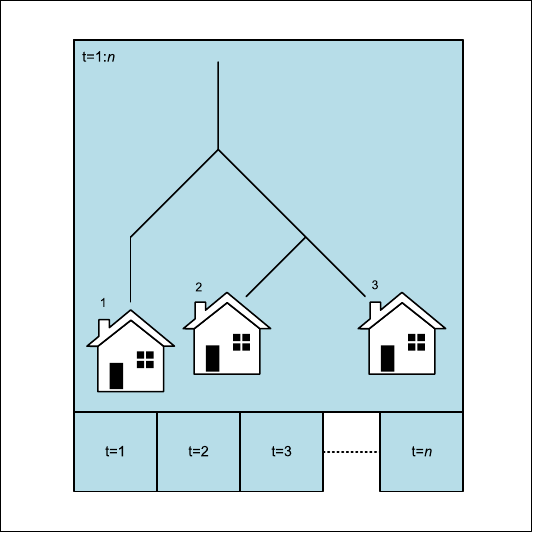}
                \caption{Central Model Structure}
                \label{fig:model_structure_central}
            \end{subfigure}
            \hfill
            \begin{subfigure}{0.63\textwidth}
                \includegraphics[width=\textwidth]{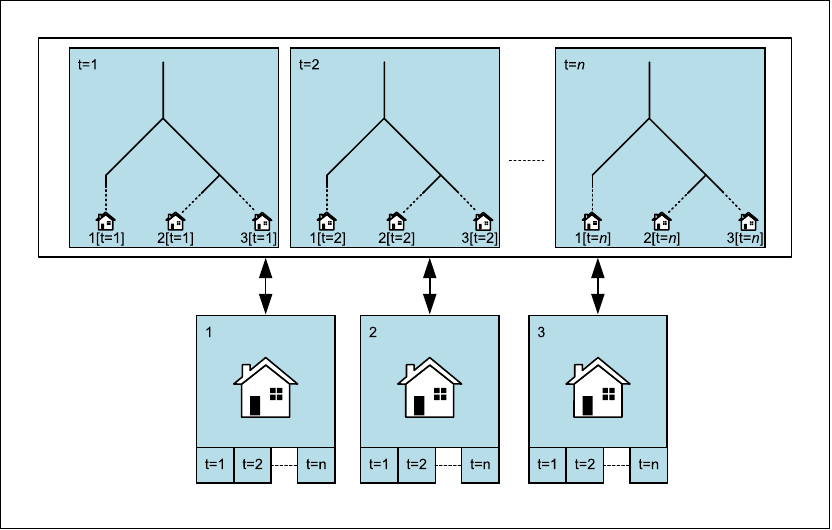}
                \caption{Decomposed Model Structure}
                \label{fig:model_structure_decomposed}
            \end{subfigure}
            \caption{Problem Decomposition Structure\label{fig:problem_decomposition_structure}}
        \end{figure}

		The full problem contains both spatial and temporal constraints. Spatial constraints govern the flow of power through the network, whilst temporal constraints link siting/sizing decisions across all timepoints and govern the charge/discharge behaviour of batteries installed within the network. It can be observed that if the network is decomposed spatially, only those partitions containing loads will be temporally constrained. No-load partitions become trivially separable in time as the solution at each timepoint is only constrained by power flow from other subproblems at a given timepoint. This property of the decomposition structure can be exploited to balance subproblem solution times, an important consideration in \gls{admm} \cite{Guo2017}, due to the aforementioned synchronicity, whereby no-load subproblems which contain significantly more buses but are only solved for a single timepoint can be balanced against subproblems containing loads which require solving across all timepoints in a single subproblem. This decomposition structure is illustrated in Figure \ref{fig:problem_decomposition_structure}, with the centralised model in Figure \ref{fig:model_structure_central} and the decomposed structure in Figure \ref{fig:model_structure_decomposed}. It should be noted that for no-load partitions, the optimisation problem becomes a feasibility problem as there are no loads present to optimise.

    \subsubsection{Partitioning Method}
\begin{figure}[H]
            \centering
            \includegraphics[width=0.55\linewidth]{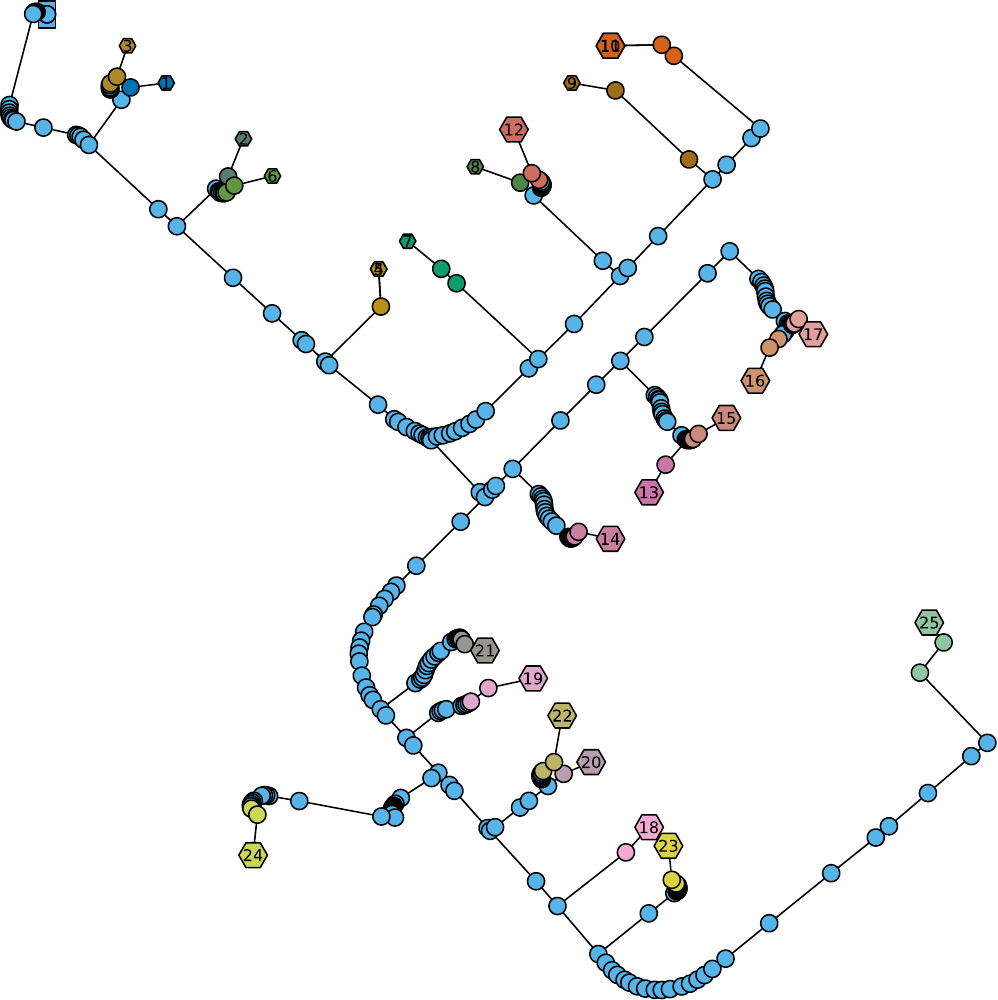}
            \caption{Partitioning of Example 25-Load Network (partitions marked in colour)}
            \label{fig:IEEE_EU_LV_25_partitioned_optBus}
        \end{figure}

        It has been shown by Guo et al. \cite{Guo2017} that the partitioning strategy employed can significantly affect the convergence performance of the \gls{admm} algorithm. The \gls{admm} implementation utilised here is synchronous, meaning that whilst the update steps in \cref{eq:admm_x_upd} can be carried out in parallel, all of them must be completed before the algorithm can move on to the next step. Therefore, it is beneficial to partition the problem such that all subproblems ideally finish solving at the same time to minimise the inefficiency of some subproblems having to wait for others to finish before the algorithm can progress. \\
        Preliminary tests were carried out to measure the solve times of both load and no-load subproblems with varying numbers of loads and buses. Here, it was  observed that the solve time for a given subproblem depends strongly on the total number of buses present in the subproblem, more so than the number of buses containing loads, which contain additional variables and constraints for the \gls{des} portion of the problem. To model an estimate of the solve time for a given partition, a regression model was fit to the results of these solve time tests. \\
        From this, a small mixed-integer quadratic problem was formulated to partition a given network structure, minimising the difference between estimated solve times by assigning buses to each subproblem. In this problem, the number of partitions was pre-defined and each load was manually assigned to one of these partitions. If more than one load was electrically close enough to another, such that placing them in separate subproblems would cause excessive partitioning, these were placed in the same partition. From this, the variables of the problem were the number of buses in each partition (integer variable) and the estimated solve time of the partition (modelled using the regression coefficients). To minimise the iterations required to drive consensus, buses could also only be added to load subproblems up to (and not including) branching buses. In this way, each load subproblem is always connected directly to the no-load subproblem made up of the main trunk buses of the network. This was also pre-defined, by specifying the maximum number of buses that could be assigned to a given partition, whilst the minimum number of buses was set to the number of load buses in each partition or the minimum number required to connect all loads in a partition if more than one is present. An additional constraint also specified that all buses must be assigned to a partition. The objective function was then to minimise the sum of the squared differences in estimated solve time for each partition, such that the optimal assignment of buses to partitions should result in subproblems having balanced solve times. \\
        An example of partitioning in this way is demonstrated in Figure \ref{fig:IEEE_EU_LV_25_partitioned_optBus} for a $25$-load network derived from the \gls{ieee} \glsxtrlong{elvtf} (see Section \ref{test_cases}).

\subsection{Test Cases}\label{test_cases}
\begin{figure}[H]
        \centering
        \includegraphics[width=12cm]{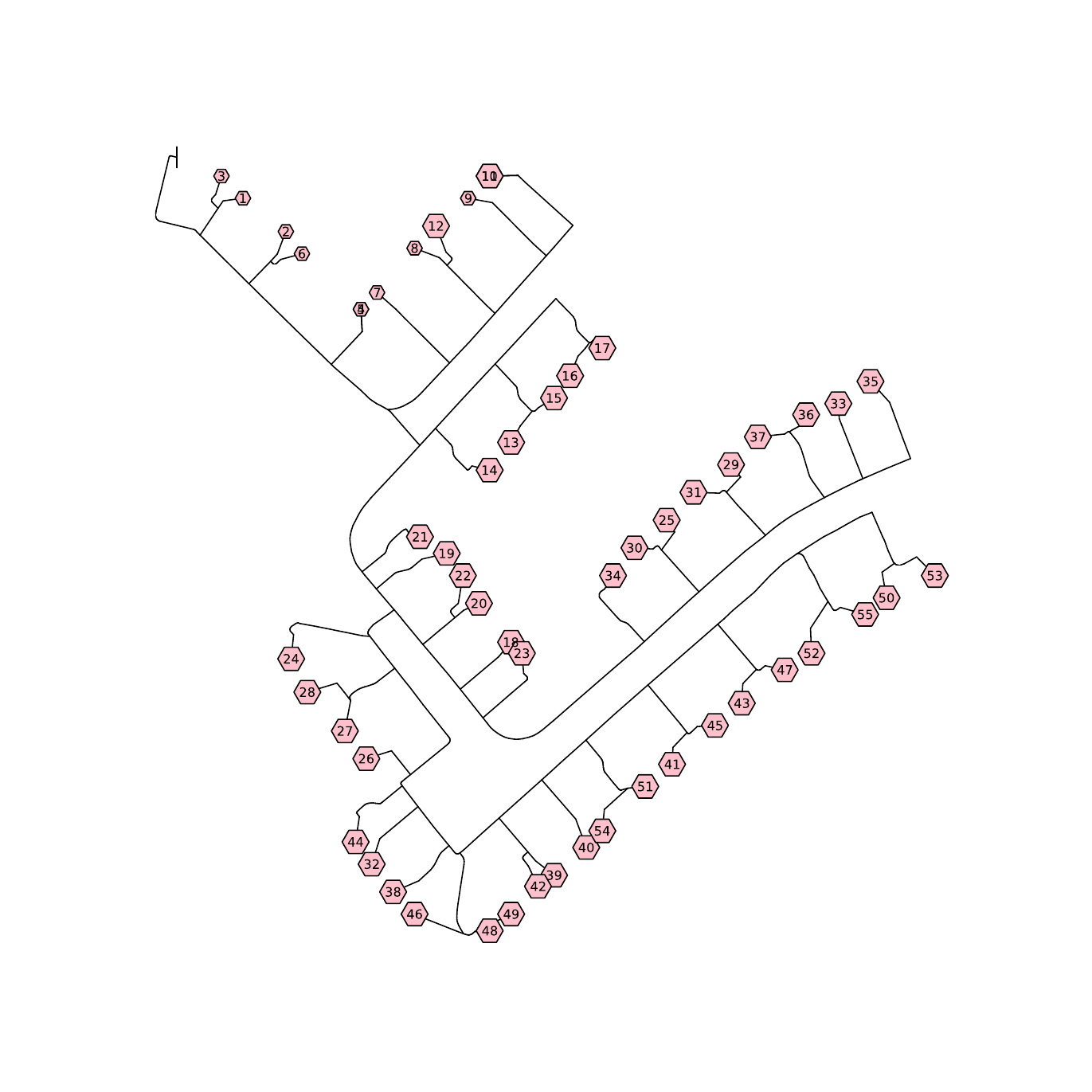}
        \caption{\Glsxtrshort{elvtf} 55-Load Test Case}
        \label{fig:elv-tf_55-load_test_case}
    \end{figure}

    The test cases considered here are derived from the \gls{ieee} \gls{elvtf} \cite{Schneider2018}. To produce sets of test feeders with varying sizes, a number of loads between $5$ and $55$ in the test feeder was specified. The network input files were then parsed such that only the subset of network elements (buses, branches, transformers etc.) required to connect this subset of loads to the feeder were kept. The corresponding number of buses and branches for the \gls{elvtf} test cases can be seen in Table \ref{tab:test_case_sizes} with the \gls{elvtf} $55$-load test network illustrated in Figure \ref{fig:elv-tf_55-load_test_case}. Note that as a result of this parsing process, some spur lines that did not lead to a load were removed from the network.

    \begin{table}[H]
        \centering
        \caption{\Glsxtrshort{elvtf} Test Case Sizes}
        \begin{tabular}{|l|l|l|}
            \hline
            \textbf{No. of Loads} & \textbf{No. of Buses} & \textbf{No. of Branches} \\
            \hline
            5 & 46 & 44 \\
            \hline
            15 & 160 & 158 \\
            \hline
            25 & 332 & 330 \\
            \hline
            35 & 456 & 454 \\
            \hline
            45 & 578 & 576 \\
            \hline
            55 & 703 & 701 \\
            \hline
        \end{tabular}
        \label{tab:test_case_sizes}
    \end{table}
	
    \subsubsection{Building Load \& Environmental Data}\label{building_load_environ}
        For the electrical and heat load input data, using a representative seasonal approach, hourly electrical and heat demand values for each building over a single $24$h period in each season are used. Electrical load data for the \gls{elvtf} test case is provided as minute-resolution data for a single winter day \cite{DeMel2024CompReform}. Using the same methodology as De Mel et al. \cite{DeMel2024CompReform}, this is averaged to provide the 24-h winter season electrical load for each building and then used to derive the electrical loads for each of the other seasons (Spring, Summer, Autumn). Heat loads were generated using Renewables.ninja \cite{Pfenninger2016,Staffell2016,RenewablesNinja} across a range of heat efficiency values between $0.1kW/^{\circ}C$ and $0.784kW/^{\circ}C$. Peak heating loads are assumed to sit within $4$kW and $9$kW, and the assumed peak heat load for each building is determined by scaling its peak electrical load between these values. This peak heat load is then matched to its closest Renewables.ninja simulation result and averaged to give the corresponding seasonal heat loads. In lieu of heat data provided by the \gls{ieee} \gls{elvtf} test case, this is intended to generate heat load data that will scale with the size of electrical load, for the purposes of selecting and sizing heating technologies. Hourly irradiance values were provided by De Mel et al. \cite{DeMel2024CompReform} and ambient temperature values were generated by Renewables.ninja and averaged into representative seasonal values \cite{Pfenninger2016,Staffell2016,RenewablesNinja}. \\
        To generate robust seasonal values, values for electrical, heat, irradiance and ambient temperature are required. For robust electrical loads, the results of Li et al. \cite{Li2018} are used, where the authors noted a standard deviation of $0.35$ for monitored mean daily electrical loads. Three standard deviations, $3\times 0.35=1.05kW$, was therefore added to the winter electrical values to give the robust seasonal electrical load. For the robust heat load, the day during winter containing the largest overall heat load was selected to be the \enquote{robust} day, representing extreme demands. The robust irradiance values are simply the average irradiance values for winter, this being the season with the lowest irradiance. For ambient temperature, the $24$h period representing the coldest day from the dataset was used, instead of the average value.

\subsection{Implementation}
    All optimisation models were formulated with the Julia programming language \cite{Bezanson2017} and JuMP software package \cite{Lubin2023}. \Gls{milp} optimisation was carried out using Gurobi (version no.: $12.0.3$) \cite{Gurobi12}, \gls{nlp} model optimisation was carried out using CONOPT 4 \cite{Drud1985} and \gls{minlp} was carried out using SBB \cite{Bussieck2001} (CONOPT 4 and SBB were provided through the GAMS solver interface version 43 \cite{GAMS2023}). Occasional solver instability was observed when running with CONOPT 4, which was attributed to bugs within the solver itself. Where this was encountered, IPOPT \cite{Wachter2006} (with the MA86 linear solver \cite{UKRI2018}) was used as a back-up solver. Some ancillary computations were carried out in Python \cite{VanRossum2009} and MATLAB \cite{MATLAB2024}. Graph analysis was carried out using Graphs.jl \cite{Graphsjl} and network plots were generated using GraphRecipes.jl \cite{Breloff2016}. Line graphs and bar charts (via Pandas \cite{Pandas1,Pandas2}) were generated using Matplotlib \cite{Hunter2007} and heatmaps were generated using Seaborn \cite{Seaborn}. Results were generated using a Dell Precision $3660$ workstation containing an Intel i9-13900K processor running at a frequency of $3.00$ GHz and with $64$ GB of installed RAM. The Microsoft Copilot generative AI tool was used to assist in code generation and debugging for the work presented in this paper, where the output was independently verified for correctness \cite{MSCopilot}. No generative AI tools were used for text in the manuscript during the writing process.

    \subsubsection{Solver Parameters}\label{solver_parameters}
        Default parameters for the Gurobi solver were used for the \gls{milp} model. For the \gls{nlp} model, the following CONOPT 4 parameters were set for the central case: \textit{Lim\_Time}: $10800$, \textit{Tol\_Scale\_Min}: $1$, \textit{Lim\_StallIter}: $100$, \textit{Tol\_Optimality}: $1\times 10^{-7}$ and for the \gls{admm} case: \textit{Lim\_Time}: $10800$, \textit{Tol\_Scale\_Min}: $1$, \textit{Lim\_StallIter}: $50$, \textit{Tol\_Optimality}: $1\times 10^{-2}$. For the \gls{minlp} model, the following SBB parameters were set: \textit{OptCR}: $0.001$, \textit{NodLim}: $5000$ and \textit{reslim}: $10800$. For IPOPT, \textit{tol}: $1\times 10^{-6}$, \textit{max\_iter}: $3000$ and \textit{max\_cpu\_time}: $10800.0$.

    \subsubsection{ADMM Parameters}
        The following \gls{admm} parameter values that were used are presented in Table \ref{tab:admm_parameter_values}, where $\beta, \zeta, \kappa$ and $\tau$ values are the same as those used by Erseghe \cite{Erseghe2015} from their \gls{ieee} feeder $123$ test case.

        \begin{table}[H]
            \centering
            \caption{\Glsxtrshort{admm} Parameter Values}
            \begin{tabular}{|l|l|}
                \hline
                \textbf{Parameter} & \textbf{Value} \\
                \hline
                $\beta$ & $10$ \\
                \hline
                $\zeta$ & $0.1$ \\
                \hline
                $\epsilon_{0}$ & $1\times 10^3$ (NLP), $5\times 10^3$ (Complementarity) \\
                \hline
                $\kappa$ & $0.99$ \\
                \hline
                $\tau$ & $1.02$ \\
                \hline
                $\lambda_{min}$ & $-1\times 10^{9}$ \\
                \hline
                $\lambda_{max}$ & $1\times 10^{9}$ \\
                \hline
                Convergence Threshold & $1\times 10^{-4}$ \\
                \hline
                Maximum Iterations & $300$ \\
                \hline
            \end{tabular}
            \label{tab:admm_parameter_values}
        \end{table}

        The convergence threshold value of $1\times 10^{-4}$ refers to the largest primal residual, defined as $\left\Vert\bm{A}\bm{x}-\bm{z}\right\Vert_{\infty}$ which is the same as the value used by Erseghe \cite{Erseghe2015} and Guo et al. \cite{Guo2017}. The issue of selecting an appropriate convergence threshold is addressed directly by Harris et al. \cite{Harris2025}, where the authors investigated the impact of changing the convergence threshold for a $500$-bus \gls{opf} test case. They observe that as the convergence threshold approached $1\times 10^{-4}$, constraint violations in the reconstructed centralised model became significant. For this reason, the primal residual convergence threshold was not raised above this value in this work.

    \section{Results}
\subsection{Centralised Model}	
\begin{figure}[H]
        \centering
        \begin{subfigure}{0.49\textwidth}
            \includegraphics[width=\textwidth]{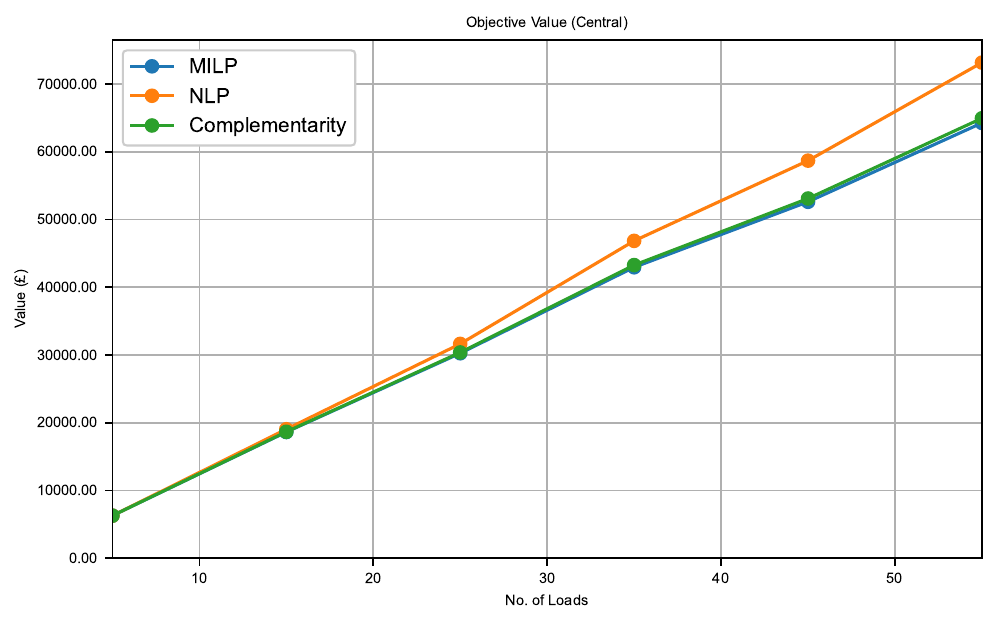}
            \caption{Central Objective Value}
            \label{fig:central_obj_val}
        \end{subfigure}
        \hfill
        \begin{subfigure}{0.49\textwidth}
            \includegraphics[width=\textwidth]{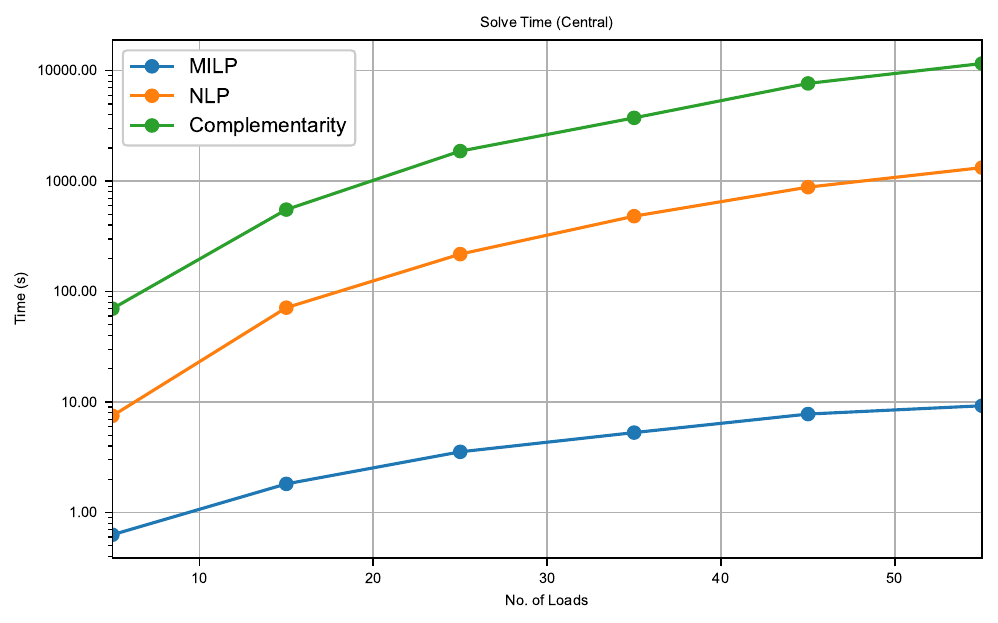}
            \caption{Central Solve Time}
            \label{fig:central_solve_time}
        \end{subfigure}
        
        \caption{Central Objective Value and Solve Time}
    \end{figure}

\begin{figure}[H]
        \centering
        \begin{subfigure}{0.49\textwidth}
            \includegraphics[width=\textwidth]{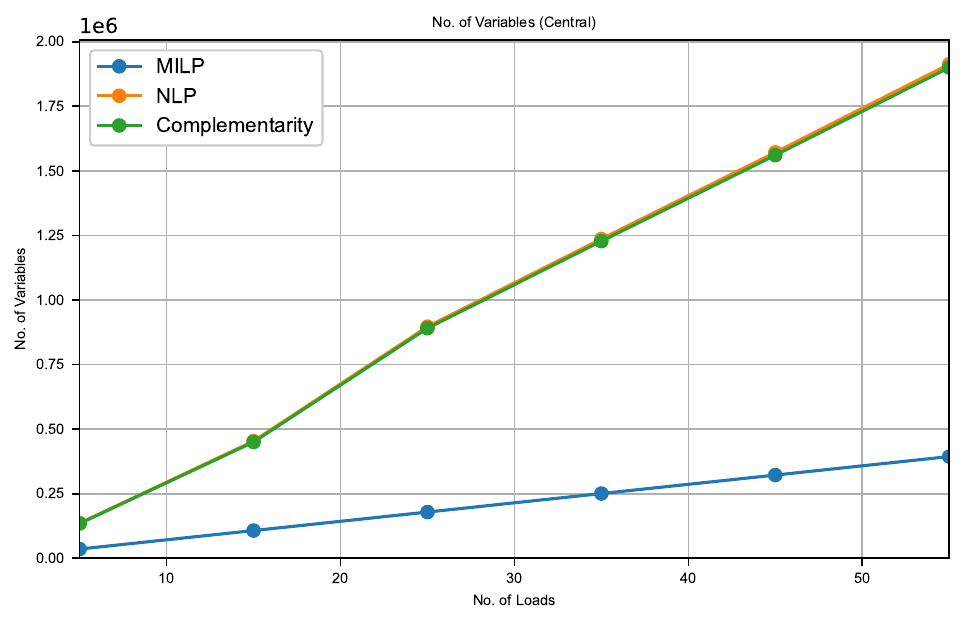}
            \caption{Central No. of Variables}
            \label{fig:central_no_of_vars}
        \end{subfigure}
        \hfill
        \begin{subfigure}{0.49\textwidth}
            \includegraphics[width=\textwidth]{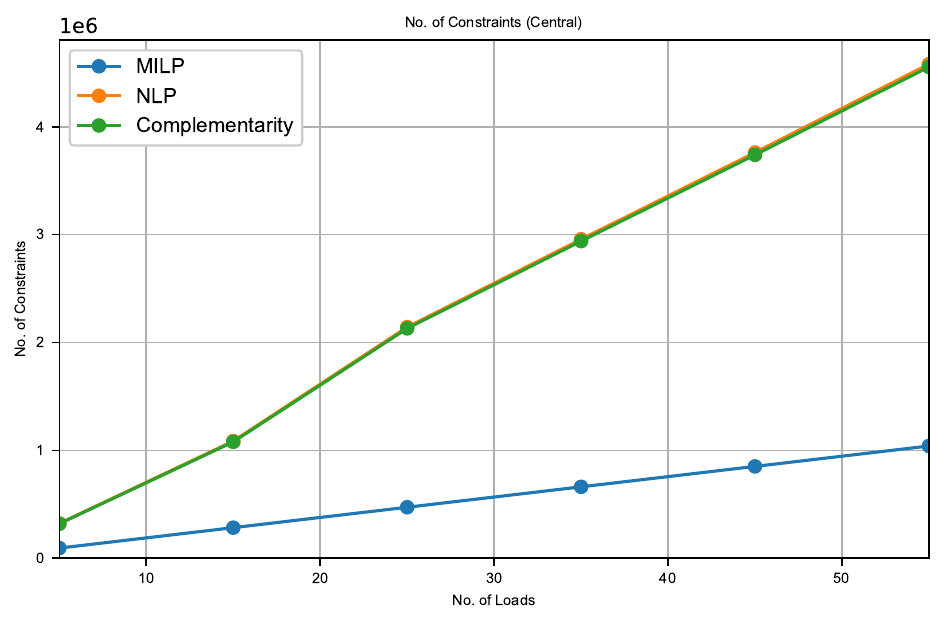}
            \caption{Central No. of Constraints}
            \label{fig:central_no_of_constrs}
        \end{subfigure}
        
        \caption{Central No. of Variables and Constraints}
    \end{figure}

    The first set of results illustrates the issue of rapidly increasing solve times, seen in Figure \ref{fig:central_solve_time} (with the corresponding objective value in Figure \ref{fig:central_obj_val}). Moving from a network containing $5$ loads to $55$ loads elicits over $2$ orders of magnitude increase in solve time, from $8$s to $1,326$s for the \gls{nlp} formulation. The increase is even more dramatic for the Complementarity formulation, moving from $70$s for $5$ loads to $11,610$s for $55$ loads. Here, the reported solve time for each formulation includes the accumulated solve time from the previous algorithmic steps (if any). Further highlighting this increase in complexity are Figures \ref{fig:central_no_of_vars} and \ref{fig:central_no_of_constrs} which illustrate the increases in the number of variables and constraints, respectively, from $134,745$ to $1,899,865$ variables and $320,600$ to $4,554,230$ constraints for the Complementarity formulation. Note that due to the implementation of the branch flow indicator constraints (\crefrange{eq:bfm_real_branch_power_flow}{eq:bfm_imag_branch_power_flow}), duplicate voltage variables are present in the model code, which have been removed from these numbers here to avoid artificially inflating the reported model size. The decision variables here are for technology installation, sizing and grid (buy/sell) and battery operation, with the number for each algorithm stage, per-load and for all timepoints, given in Table \ref{tab:model_decision_variables}. The \gls{minlp} formulation was also run, but reached the timeout without solving for all network sizes, and therefore is not displayed here.

    \begin{table}[H]
        \centering
        \caption{Model Decision Variables (per-load, $120$ timepoints)}
        \begin{tabular}{|p{2.8cm}|p{2.5cm}|p{1cm}|p{1cm}|p{3cm}|p{2.4cm}|}
        \hline
            \textbf{Variable Type} & \textbf{Component} & \textbf{\Glsxtrshort{milp}} & \textbf{\Glsxtrshort{nlp}} & \textbf{Complementarity} & \textbf{Fixed Complementarity} \\
            \hline
            \multirow{5}{*}{Technology Choice} & \Glsxtrshort{pv} & $0$ & $0$ & $0$ & $0$\\ \cline{2-6}
                                               & Battery & $1$ & $0$ & $0$ & $0$ \\ \cline{2-6}
                                               & Heat Pump & $6$ & $0$ & $0$ & $0$\\ \cline{2-6}
                                               & Hot Water Tank & $4$ & $0$ & $0$ & $0$\\ \cline{2-6}
                                               & Boiler & $1$ & $0$ & $0$ & $0$ \\ 
            \hline
            \multirow{5}{*}{Technology Size}   & \Glsxtrshort{pv} & $1$ & $1$ & $1$ & $1$\\ \cline{2-6}
                                               & Battery & $1$ & $1$ & $1$ & $1$ \\ \cline{2-6}
                                               & Heat Pump & $0$ & $0$ & $0$ & $0$\\ \cline{2-6}
                                               & Hot Water Tank & $0$ & $0$ & $0$ & $0$\\ \cline{2-6}
                                               & Boiler & $1$ & $1$ & $1$ & $1$ \\ 
            \hline
            \multirow{2}{*}{Operation}         & Battery & $240$ & $120$ & $240$ & $120$ \\ \cline{2-6}
                                               & Grid & $240$ & $120$ & $240$ & $120$ \\
            \hline
        \end{tabular}
        \label{tab:model_decision_variables}
    \end{table}

\begin{figure}[H]
        \centering
        \includegraphics[scale=0.65]{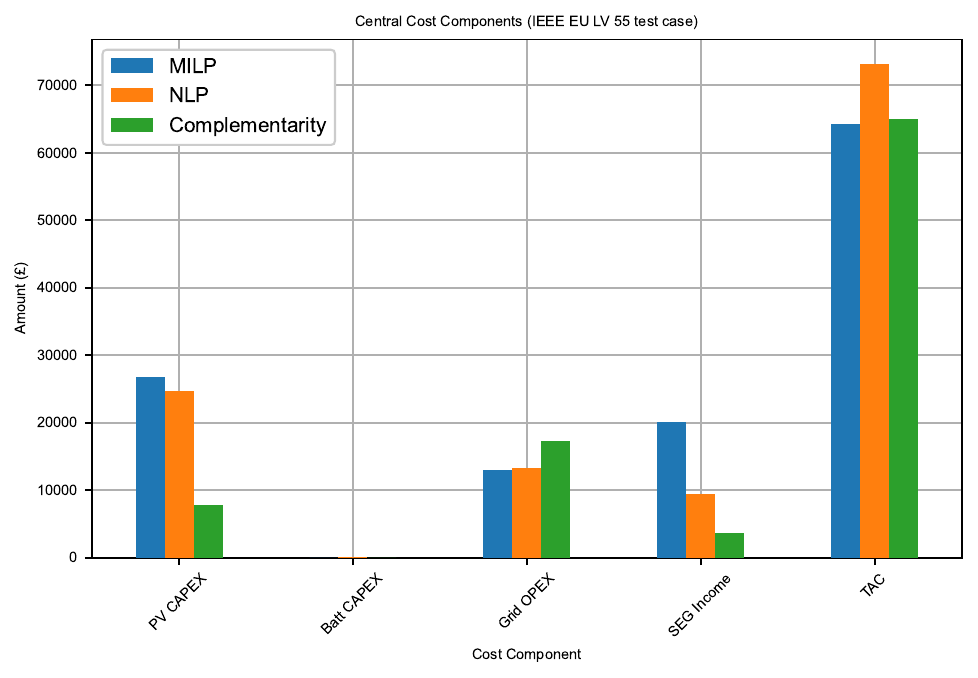}
        \caption{Central 55 Load Cost Components}
        \label{fig:central_55_load_cost_comps}
    \end{figure}

    Also of note from these results is the increasing disparity between the \gls{milp}, \gls{nlp} and Complementarity solutions. As expected, the \gls{milp} solution (step \ref{line:central_milp_solve} in Algorithm \ref{alg:central_solve}) returns effectively a lower bound as it is not constrained by power flow in the network. Adding these constraints in the \gls{nlp} formulation (step \ref{line:central_nlp_solve} in Algorithm \ref{alg:central_solve}) results in a high objective value, which can then be lowered almost to the same value again in the Complementarity formulation (step \ref{line:central_comp_solve} in Algorithm \ref{alg:central_solve}) where the model has additional degrees of freedom to adjust for these power flow constraints. \\
    Examining the cost components in more detail reveals the mechanism for this. In Figure \ref{fig:central_55_load_cost_comps}, a relevant subsection of the model costs are shown. Moving from \gls{milp}$\rightarrow$\gls{nlp}$\rightarrow$Complementarity shows large reductions in \gls{pv} installation, with corresponding decreases in \gls{seg} income. In the \gls{milp} model, each load has the unfettered ability to export power at peak times, which comes with the corresponding income. This is reduced in the \gls{nlp} formulation, when network limit constraints may prevent loads from all simultaneously exporting in the middle of the day. This is further reduced in the Complementarity formulation where the model is able to adjust the operational schedules in the model, allowing it to meet the loads with a drastically reduced amount of installed \gls{pv} capacity. Note: for brevity in the figure the fixed operating costs for \gls{pv} and battery were omitted, alongside boiler, heat pump and hot water tank costs which remain constant across formulations. In this case, across all network sizes the solution opted to install negligible battery storage and gas boilers instead of heat pumps/hot water tanks. Carbon emissions were not considered in this model, which if they were, may have altered the decision by the model to install boilers over heat pumps.

\begin{figure}[H]
        \centering
        \begin{subfigure}{0.49\textwidth}
            \includegraphics[width=\textwidth]{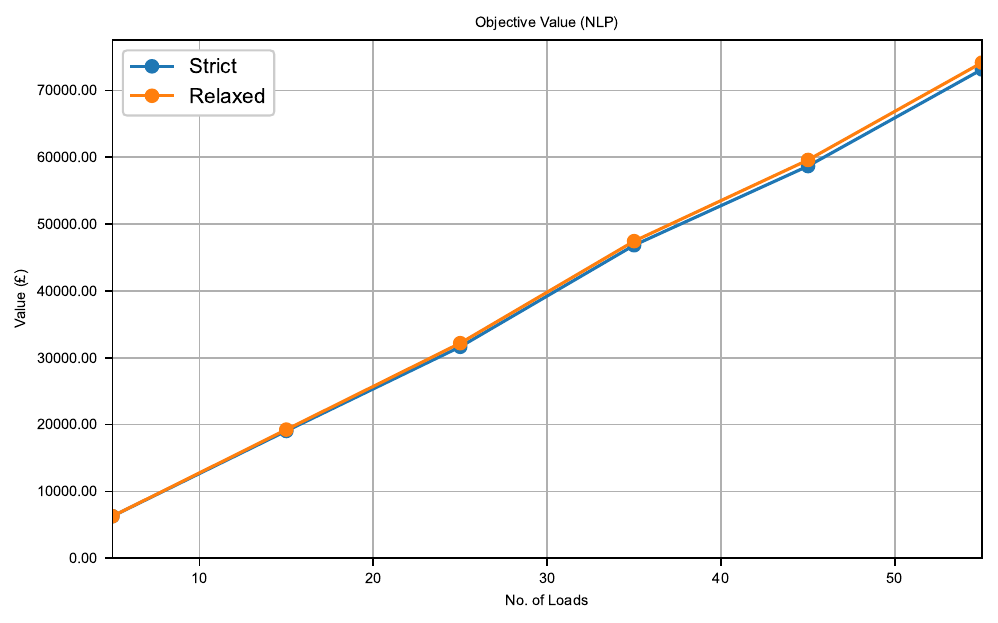}
            \caption{Objective Value}
            \label{fig:central_nlp_strict_vs_relaxed_obj_val}
        \end{subfigure}
        \hfill
        \begin{subfigure}{0.49\textwidth}
            \includegraphics[width=\textwidth]{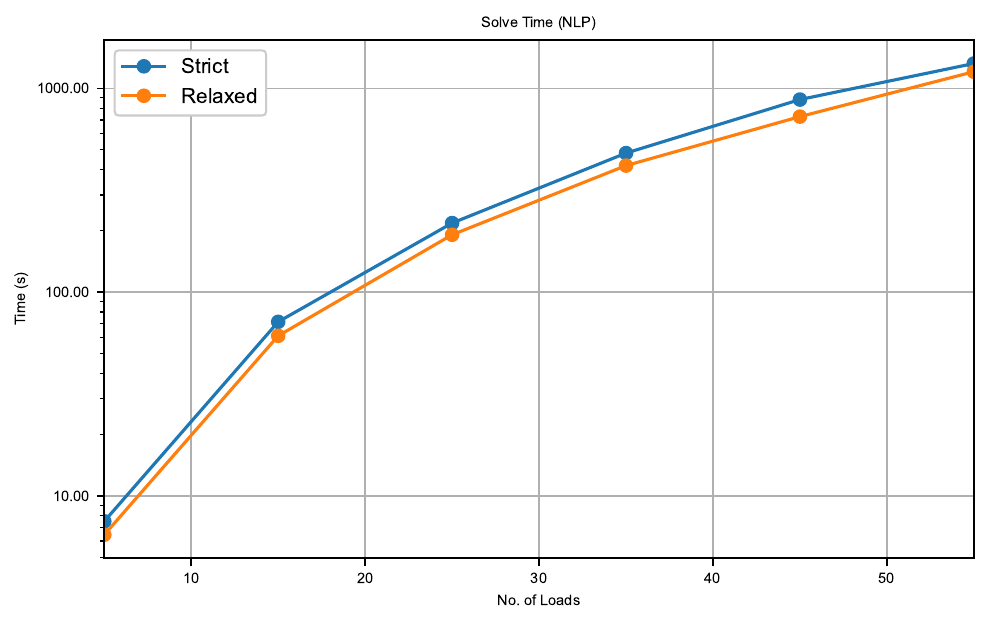}
            \caption{Solve Time}
            \label{fig:central_nlp_strict_vs_relaxed_solve_time}
        \end{subfigure}
        
        \caption{Central \glsxtrshort{nlp} Solve Strict vs Relaxed Solver Parameters}
    \end{figure}

\begin{figure}[H]
        \centering
        \begin{subfigure}{0.49\textwidth}
            \includegraphics[width=\textwidth]{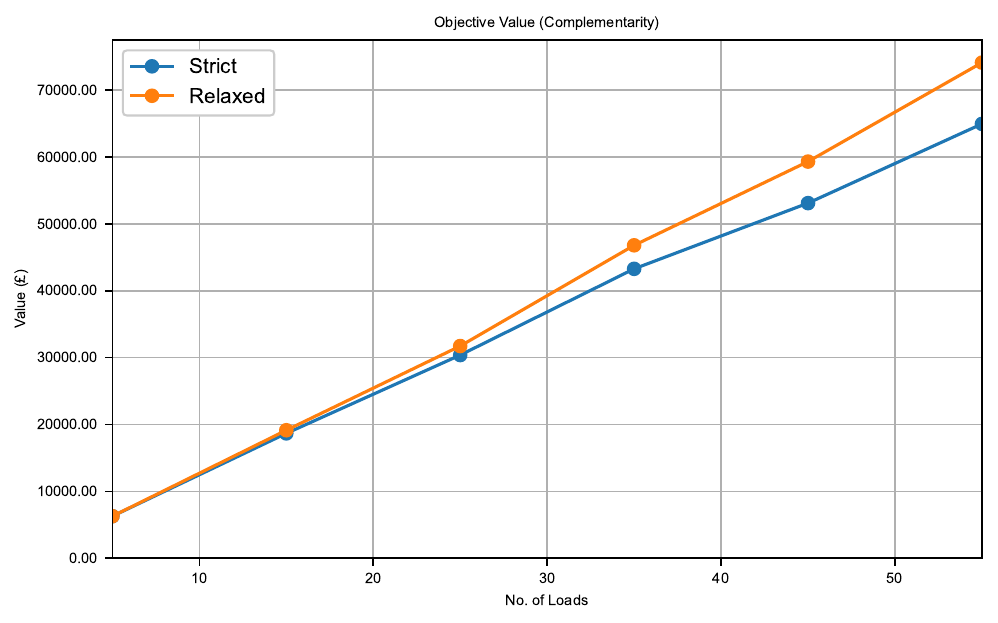}
            \caption{Objective Value}
            \label{fig:central_comp_strict_vs_relaxed_obj_val}
        \end{subfigure}
        \hfill
        \begin{subfigure}{0.49\textwidth}
            \includegraphics[width=\textwidth]{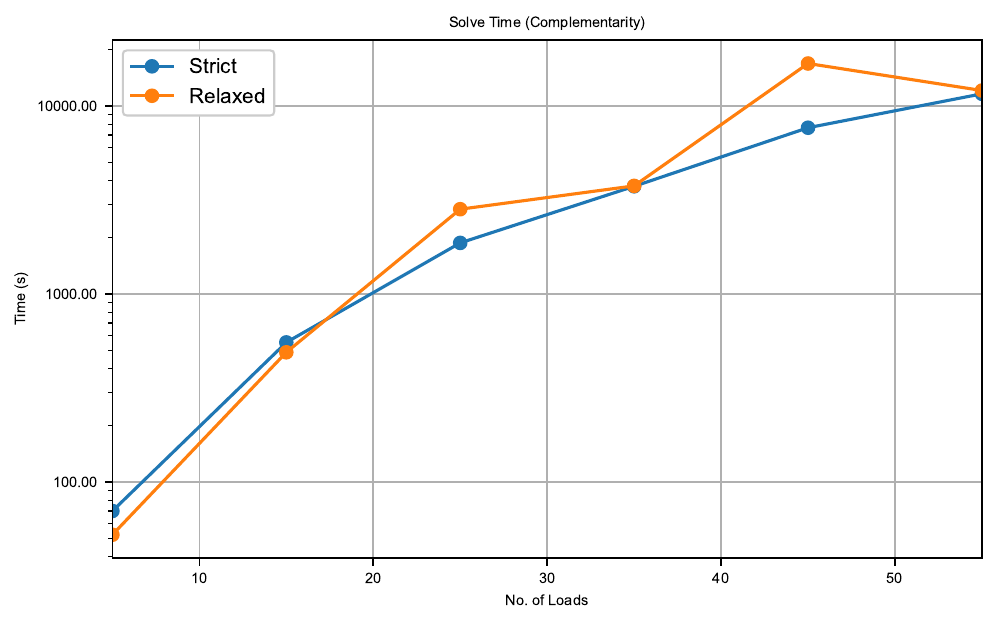}
            \caption{Solve Time}
            \label{fig:central_comp_strict_vs_relaxed_solve_time}
        \end{subfigure}
        
        \caption{Central Complementarity Solve Strict vs Relaxed Solver Parameters}
    \end{figure}

    As is noted in Section \ref{solver_parameters}, relaxed solver parameters were used when solving using \gls{admm} vs centrally. To determine the effect on solver performance that these would have had if run centrally, the problem was solved centrally using both strict and relaxed parameters. The comparison in terms of objective value is shown in Figure \ref{fig:central_nlp_strict_vs_relaxed_obj_val} and solve time in Figure \ref{fig:central_nlp_strict_vs_relaxed_solve_time} for the \gls{nlp} formulation and Figures \ref{fig:central_comp_strict_vs_relaxed_obj_val} and \ref{fig:central_comp_strict_vs_relaxed_solve_time} respectively for the Complementarity formulation. Relaxing the solver optimality parameters for the \gls{nlp} formulation resulted in an optimality gap of $1.38\%$ for $55$ loads, with a corresponding solve time ratio of $0.91$. For the Complementarity formulation, the optimality gap was $14.10\%$ with a solve time ratio of $1.04$.

\subsection{ADMM Model}
    \subsubsection{Penalty Parameter Initial Value Selection}
\begin{figure}[H]
            \centering
            \includegraphics[scale=0.7]{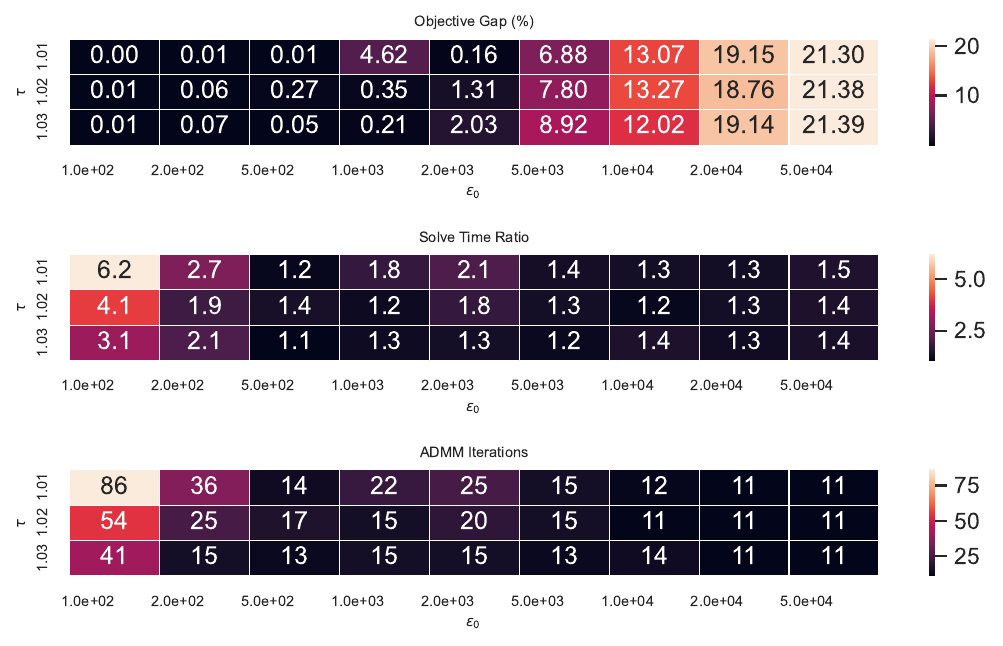}
            \caption{\Glsxtrshort{elvtf} \glsxtrshort{nlp} 5-Load Parameter Sweep}
            \label{fig:admm_param_sweep_nlp_5_load}
        \end{figure}

\begin{figure}[H]
            \centering
            \includegraphics[scale=0.7]{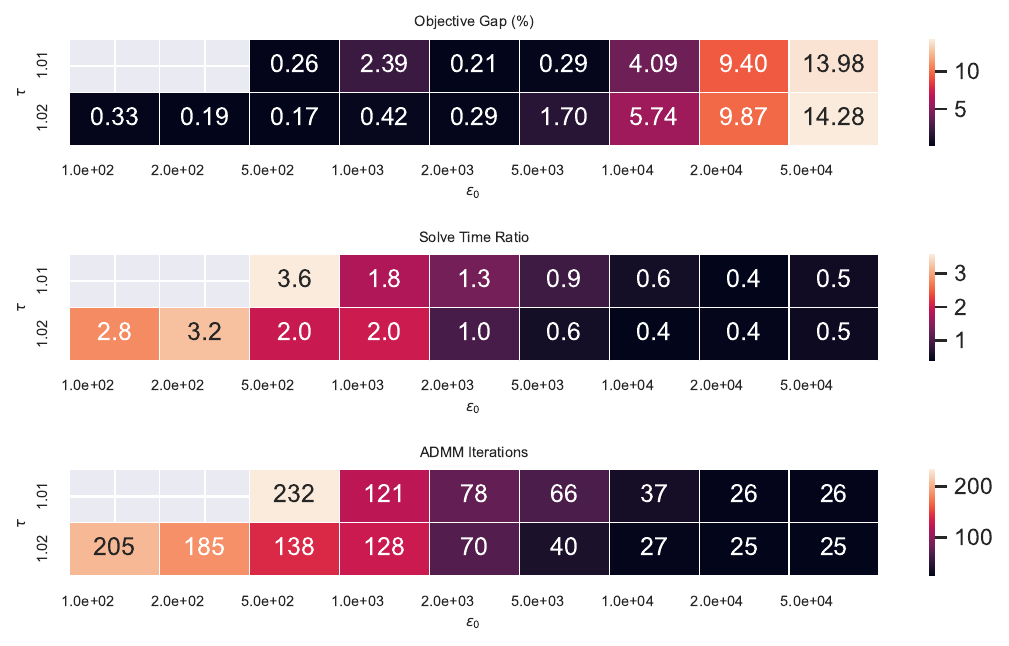}
            \caption{\Glsxtrshort{elvtf} \glsxtrshort{nlp} 15-Load Parameter Sweep}
            \label{fig:admm_param_sweep_nlp_15_load}
        \end{figure}

\begin{figure}[H]
            \centering
            \includegraphics[scale=0.7]{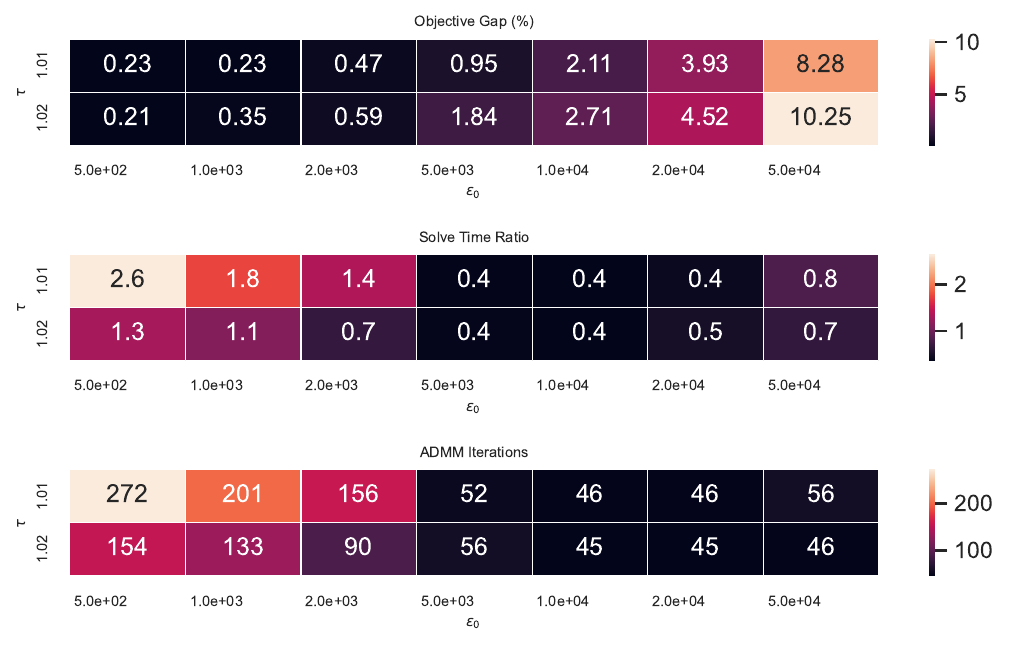}
            \caption{\Glsxtrshort{elvtf} \glsxtrshort{nlp} 25-Load Parameter Sweep}
            \label{fig:admm_param_sweep_nlp_25_load}
        \end{figure}

\begin{figure}[H]
            \centering
            \includegraphics[scale=0.7]{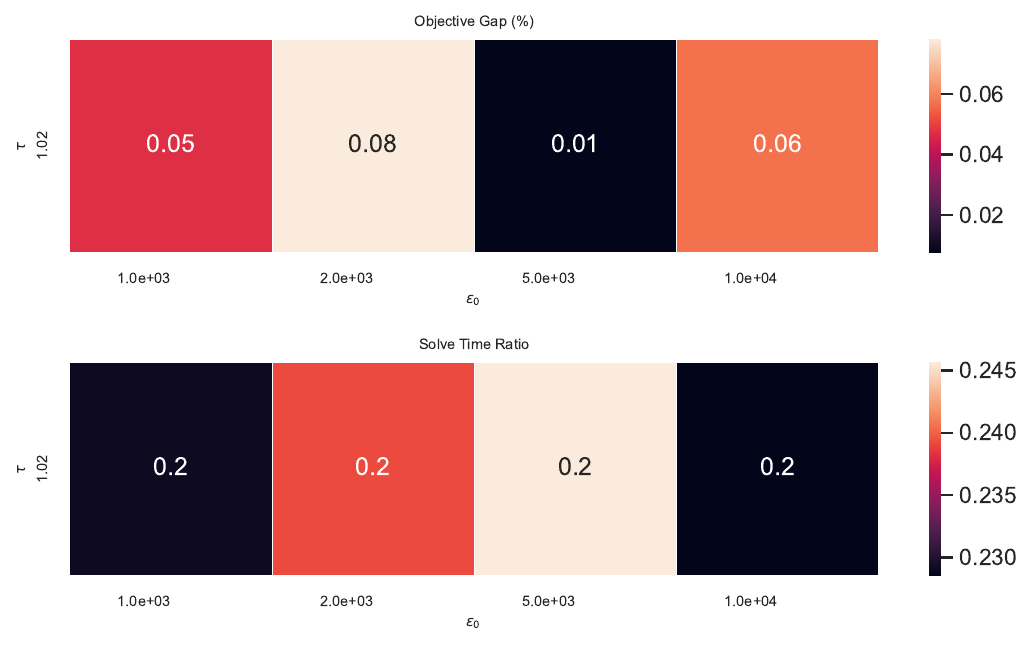}
            \caption{\Glsxtrshort{elvtf} Complementarity 5-Load Parameter Sweep}
            \label{fig:admm_param_sweep_comp_5_load}
        \end{figure}

\begin{figure}[H]
            \centering
            \includegraphics[scale=0.7]{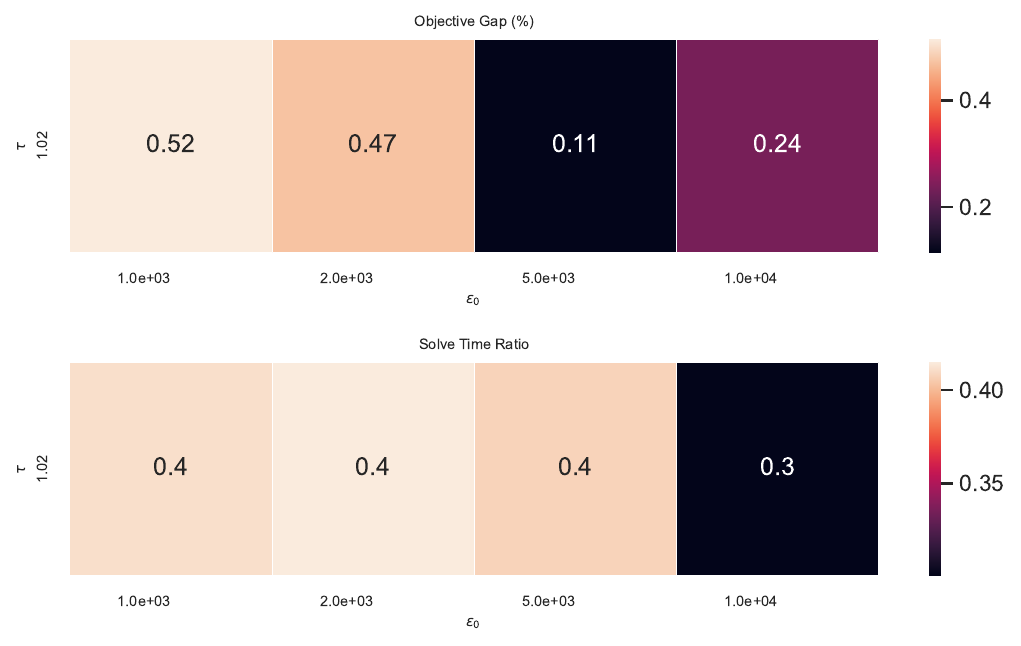}
            \caption{\Gls{elvtf} Complementarity 15-Load Parameter Sweep}
            \label{fig:admm_param_sweep_comp_15_load}
        \end{figure}

\begin{figure}[H]
            \centering
            \includegraphics[scale=0.7]{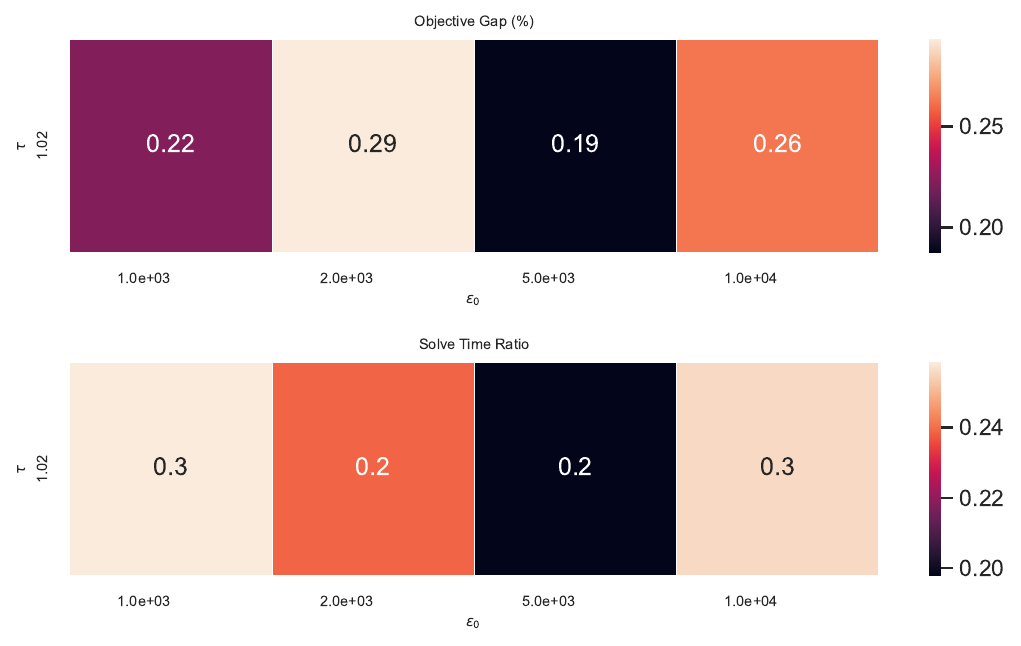}
            \caption{\Glsxtrshort{elvtf} Complementarity 25-Load Parameter Sweep}
            \label{fig:admm_param_sweep_comp_25_load}
        \end{figure}

        As noted by Guo et al. \cite{Guo2017}, optimal parameter selection for \gls{admm} is normally carried out through experimentation. It is noted that the \gls{admm} penalty parameter (in this case, the initial penalty parameter $\epsilon_{0}$) has a high impact on how the algorithm performs, giving preference between solution quality and convergence speed. By examining \cref{eq:admm_aug_lagrange}, it can be seen that a higher penalty parameter results in larger values for the augmented term, encouraging solution feasibility at the expense of optimality. Conversely, a lower penalty parameter value prioritises the initial objective terms, potentially leading to a better quality solution but taking longer to reach as solution feasibility is no longer encouraged to the same degree. \\
        The aim therefore is to find a value for $\epsilon_{0}$ that is satisfactory for both performance metrics. This was done using a parameter sweep for $5$, $15$ and $25$ loads, with the results shown in Figures \ref{fig:admm_param_sweep_nlp_5_load}, \ref{fig:admm_param_sweep_nlp_15_load} and \ref{fig:admm_param_sweep_nlp_25_load} respectively for the \gls{nlp} formulation. Due to the long experiment run times for the larger networks and $\epsilon_{0}$ values, the swept values were updated between network sizes. The empty entries in Figure \ref{fig:admm_param_sweep_nlp_15_load} indicate that the \gls{admm} algorithm did not converge within the specified number of iterations. From these results, the value of $\epsilon_{0}=1\times 10^{3}$ was selected as a reasonable trade-off between objective gap and solve time ratio, when compared against the centralised formulation (this compares against values of $10-30$ used by Erseghe \cite{Erseghe2015}, $1\times 10^4-1\times 10^7$ used by Guo et al. \cite{Guo2017} and $1$ by Nick et al. \cite{Nick2015}). \\
        Following this, $\epsilon_{0}$ was fixed for the \gls{nlp} formulation and a further sweep was carried out for the Complementarity formulation (here carried out with no fixed complementarity step), with results shown in Figures \ref{fig:admm_param_sweep_comp_5_load}, \ref{fig:admm_param_sweep_comp_15_load} and \ref{fig:admm_param_sweep_comp_25_load}. These indicated that $\epsilon_{0}$ could be set higher, reducing the solve time ratio, whilst still achieving an acceptable objective gap. From these results, a value of $\epsilon_{0}=5\times 10^{3}$ was selected. \\
        As can be seen from the \gls{nlp} sweep figures, different values of the $\tau$ parameter, controlling how quickly the penalty parameter is increased as the algorithm progresses \cite{Erseghe2015}, was also tested. From the tested results, it was decided to keep the default value of $1.02$. \\
        Here, objective gap is defined as:
        \begin{equation}\label{eq:objective_gap}
            obj\_gap=\left(\frac{obj_{ADMM}-obj_{cent}}{obj_{cent}}\right)\times 100
        \end{equation}
        where $obj_{ADMM}$ is the objective value for the \gls{admm} solution and $obj_{cent}$ is the objective value for the centralised solution. The solve time ratio is defined as:
        \begin{equation}\label{eq:solve_time}
            t_{ratio}=\frac{t_{ADMM}}{t_{cent}}
        \end{equation}
        where $t_{ADMM}$ is the \gls{admm} solution solve time and $t_{cent}$ is the centralised solution solve time. The solution solve time for \gls{admm}, $t_{ADMM}$, is the sum of all \gls{admm} iteration solve times, each of which is the maximum solver time taken over all subproblems (assumed to be solved simultaneously in parallel). The time taken to communicate results between subproblems is not considered here, in line with Erseghe \cite{Erseghe2015} and Guo et al. \cite{Guo2017}, as all computation is carried out locally on a single machine. As well as this, the time taken for the $z$ and $\lambda$ updates (\cref{eq:admm_z_upd} and \cref{eq:admm_lambda_upd} respectively), which mainly constitute sparse matrix-vector multiplications, was taken to be negligible here and not considered as the time taken for these multiplications was benchmarked to be less than $5.5$ms, compared to subproblem solution times of the order of seconds for larger networks.
        
    \subsubsection{ADMM Results}
\begin{figure}[H]
            \centering
            \begin{subfigure}{0.49\textwidth}
                \includegraphics[width=\textwidth]{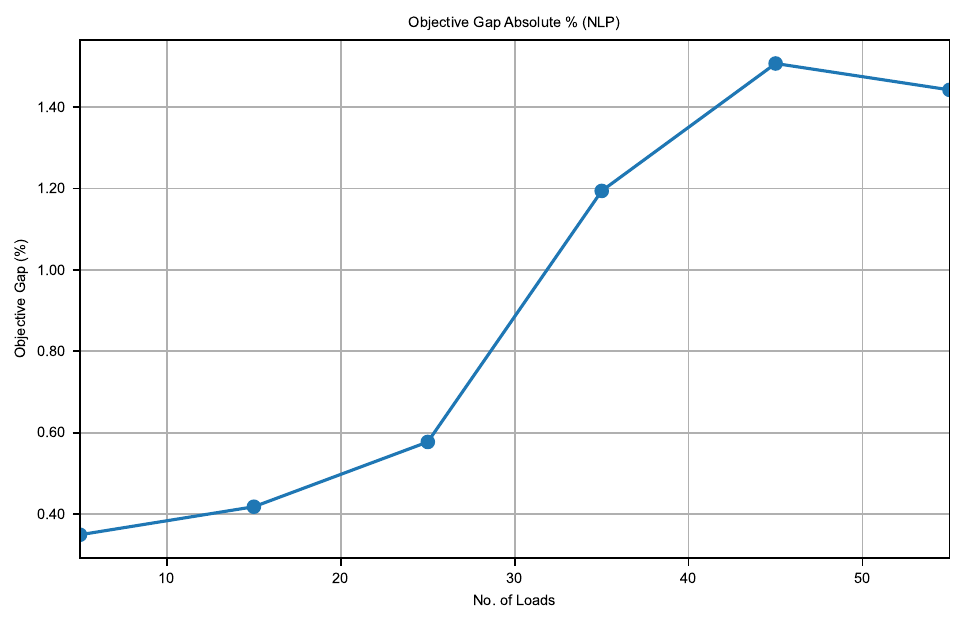}
                \caption{Objective Gap Percentage}
                \label{fig:dist_NLP_obj_gap}
            \end{subfigure}
            \hfill
            \begin{subfigure}{0.49\textwidth}
                \includegraphics[width=\textwidth]{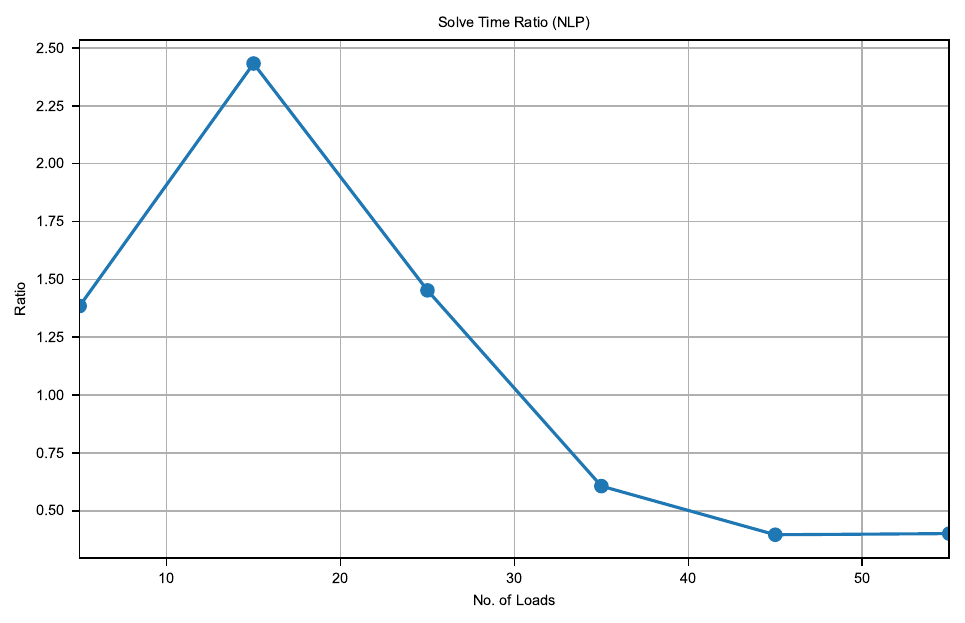}
                \caption{Solve Time Ratio}
                \label{fig:dist_NLP_solve_time_ratio}
            \end{subfigure}
            
            \caption{\Glsxtrshort{admm} \glsxtrshort{nlp} Objective Gap and Solve Time Ratio}
        \end{figure}

\begin{figure}[H]
            \centering
            \begin{subfigure}{0.49\textwidth}
                \includegraphics[width=\textwidth]{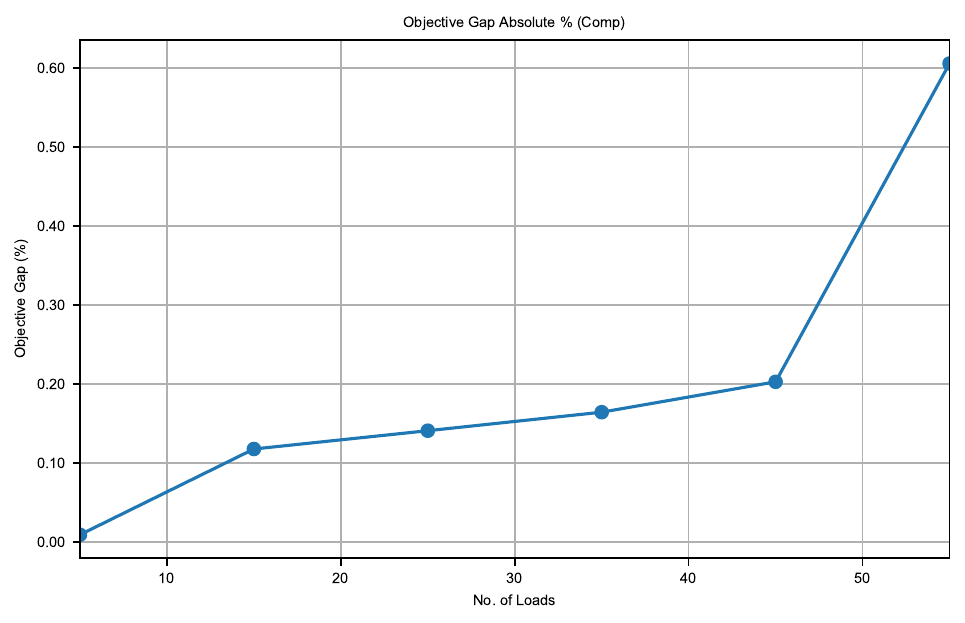}
                \caption{Objective Gap}
                \label{fig:dist_comp_obj_gap}
            \end{subfigure}
            \hfill
            \begin{subfigure}{0.49\textwidth}
                \includegraphics[width=\textwidth]{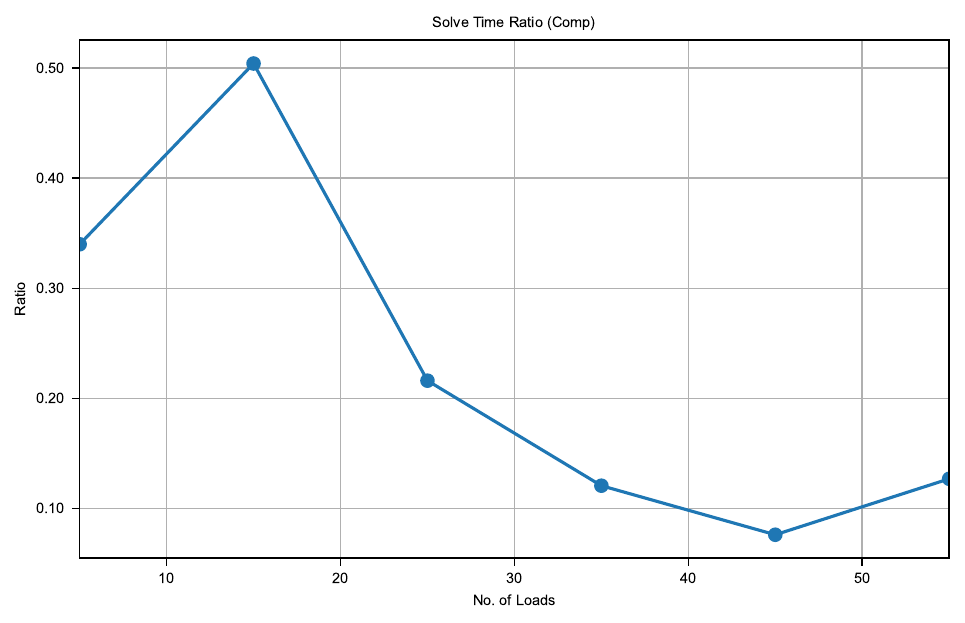}
                \caption{Solve Time Ratio}
                \label{fig:dist_comp_solve_time_ratio}
            \end{subfigure}
            
            \caption{\Glsxtrshort{admm} Complementarity Objective Gap Percentage and Solve Time Ratio}
        \end{figure}

\begin{figure}[H]
            \centering
            \begin{subfigure}{0.49\textwidth}
                \includegraphics[width=\textwidth]{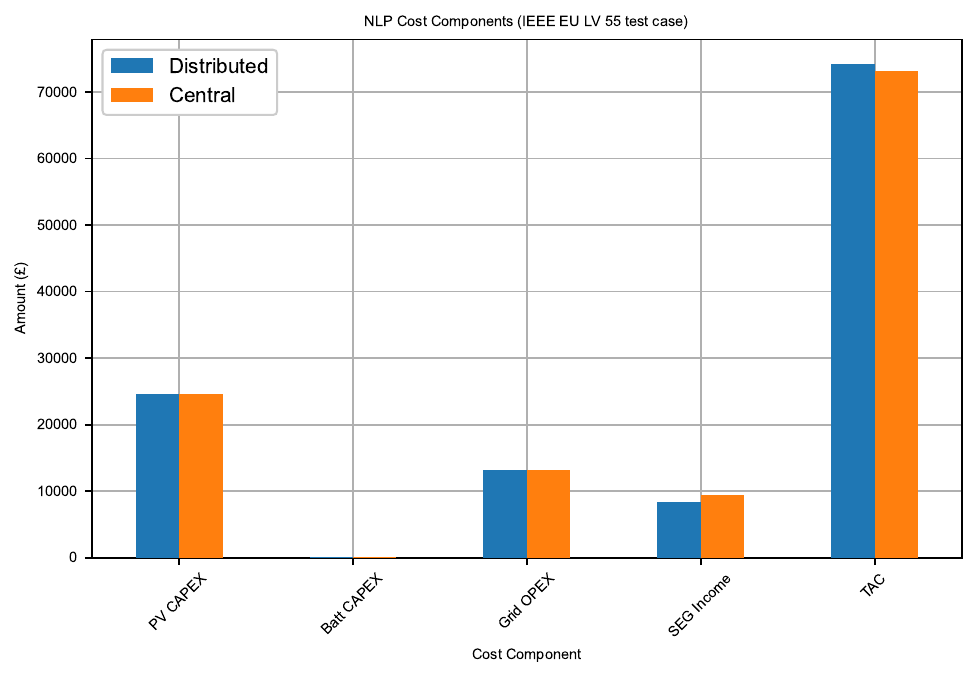}
                \caption{\Glsxtrshort{nlp} Cost Comparison}
                \label{fig:55_load_NLP_dist_central_cost_comparison}
            \end{subfigure}
            \hfill
            \begin{subfigure}{0.49\textwidth}
                \includegraphics[width=\textwidth]{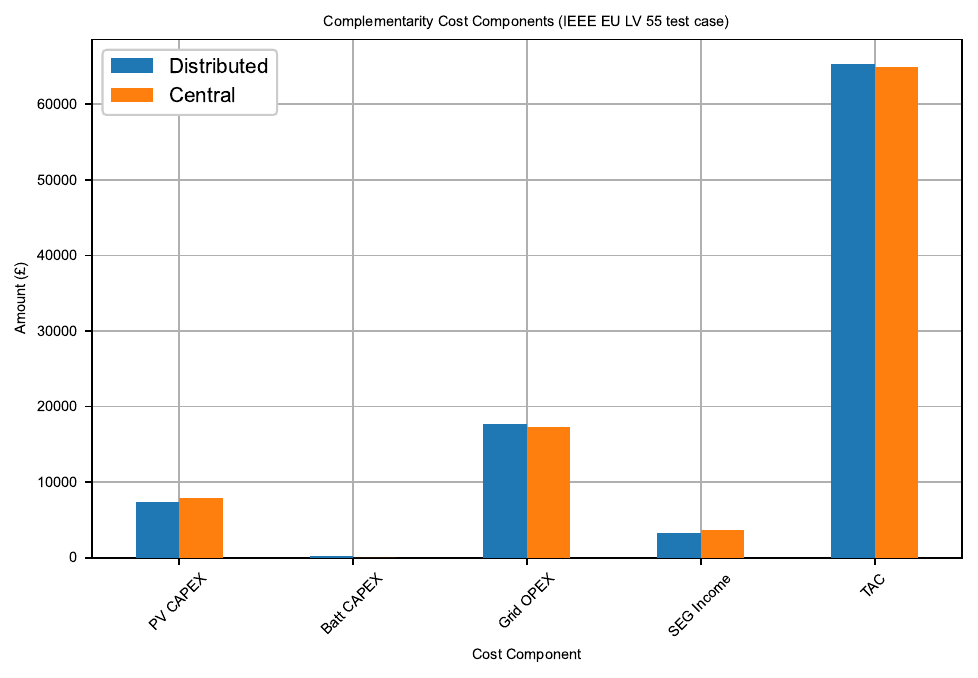}
                \caption{Complementarity Cost Comparison}
                \label{fig:55_load_comp_dist_central_cost_comparison}
            \end{subfigure}
            
            \caption{\Glsxtrshort{elvtf} 55-Load \glsxtrshort{nlp} and Complementarity Centralised vs Distributed Cost Component Comparison}
        \end{figure}

        The main results for this paper can be seen in Figures \ref{fig:dist_NLP_obj_gap}, \ref{fig:dist_NLP_solve_time_ratio}, \ref{fig:dist_comp_obj_gap} and \ref{fig:dist_comp_solve_time_ratio}. Here, the \gls{admm} approach is tested against the centralised approach (the minimum of the values for the strict and relaxed central solves) in terms of both optimality gap (Figures \ref{fig:dist_NLP_obj_gap} and \ref{fig:dist_comp_obj_gap}) and solve time ratio (Figures \ref{fig:dist_NLP_solve_time_ratio} and \ref{fig:dist_comp_solve_time_ratio}), for both the \gls{nlp} and Complementarity formulations. \\
        For both the \gls{nlp} and Complementarity solves, the solve time ratio drastically reduces as the size of the network grows, reaching $0.40$ for \gls{nlp} and $0.13$ for Complementarity for the network containing $55$ loads. Whilst the objective gap for \gls{nlp} increases for larger network sizes, reaching a maximum of $1.51\%$ for $45$ loads, for most networks this is almost entirely mitigated when moving to the Complementarity formulation, with a maximum gap of $0.61\%$ seen across the tested networks. The differences in solution cost components for distributed and centralised are elucidated in Figures \ref{fig:55_load_NLP_dist_central_cost_comparison} and \ref{fig:55_load_comp_dist_central_cost_comparison} for \gls{nlp} and Complementarity respectively. \\
        The results shown here highlight the ability of the model to provide optimisation-based studies and decision-support for \gls{des} design problems involving larger networks with unbalanced three-phase power flow, whilst avoiding the use of power flow approximations. In doing so, it enables a number of different aspects of the model to be studied. These include studying the impact of the installation of \glspl{der} in a network in terms of both capital and operating costs, whether the network structure itself impacts these \gls{der} installations, for example if voltage limits are likely to be reached during operation, as well as the impact of changing model/network parameters. Running in a distributed manner using \gls{admm} improves the scalability of this when studying larger networks, through the use of the hybrid spatial-temporal decomposition of the design problem space.

    \subsubsection{Power Flow Validation}
        Spatially decomposing the network leads to the possibility of constraint violations (see Harris et al. \cite{Harris2025}) due to the relaxation of bus variable equality across tie-line branches. Therefore, the results obtained by the distributed solutions were validated by carrying out power flow calculations in Pandapower \cite{Thurner2018} (with Python \cite{VanRossum2009} via PyCall.jl \cite{PyCall}) using the distributed solution real and reactive power injections. The power flow solution was then analysed for any voltage magnitude constraint violations (voltage angle and line flow limits were not considered here). Statistics pertaining to the voltage magnitude upper constraint violation are presented for the \gls{nlp} and Complementarity steps in Tables \ref{tab:nlp_v_mag_upper_constraint_viol} and \ref{tab:comp_v_mag_upper_constraint_viol} respectively. Across the test cases, the maximum observed constraint violation was $0.205$\% for the \gls{nlp} step but only $0.086$\% for the Complementarity step, which was deemed to be acceptable. Voltage magnitude lower constraint violations were also analysed but no constraint violations were found across any test case. Note that average constraint violation statistics include buses with no constraint violations.

        \begin{table}[H]
            \centering
            \caption{\Glsxtrshort{nlp} Voltage Magnitude Upper Limit Constraint Violation Statistics}
            \begin{tabular}{|p{1.6cm}|p{2.2cm}|p{1.9cm}|p{2.8cm}|}
                \hline
                \textbf{No. of Loads} & \textbf{Average Constraint Violation (\%)} & \textbf{Max. Constraint Violation (\%)} & \textbf{Percentage of Constraints Violated (\%)} \\
                \hline
                $5$ & $0.000000$ & $0.00$ & $0.00$ \\
                \hline
                $15$ & $0.000093$ & $0.051$ & $0.53$ \\
                \hline
                $25$ & $0.000181$ & $0.070$ & $1.14$ \\
                \hline
                $35$ & $0.000983$ & $0.177$ & $2.89$ \\
                \hline
                $45$ & $0.001710$ & $0.167$ & $4.59$ \\
                \hline
                $55$ &$ 0.004164$ & $0.205$ & $7.84$ \\
                \hline
            \end{tabular}
            \label{tab:nlp_v_mag_upper_constraint_viol}
        \end{table}

        \begin{table}[H]
            \centering
            \caption{Complementarity Voltage Magnitude Upper Limit Constraint Violation Statistics}
            \begin{tabular}{|p{1.6cm}|p{2.2cm}|p{1.9cm}|p{2.8cm}|}
                \hline
                \textbf{No. of Loads} & \textbf{Average Constraint Violation (\%)} & \textbf{Max. Constraint Violation (\%)} & \textbf{Percentage of Constraints Violated (\%)} \\
                \hline
                $5$ & $0.0000000$ & $0.00$ & $0.00$ \\
                \hline
                $15$ & $0.0000021$ & $0.022$ & $0.023$ \\
                \hline
                $25$ & $0.0000280$ & $0.081$ & $0.124$ \\
                \hline
                $35$ & $0.0000473$ & $0.067$ & $0.201$ \\
                \hline
                $45$ & $0.0000770$ & $0.071$ & $0.514$ \\
                \hline
                $55$ & $0.0002431$ & $0.086$ & $1.024$ \\
                \hline
            \end{tabular}
            \label{tab:comp_v_mag_upper_constraint_viol}
        \end{table}

    \section{Conclusion}
Introduction of large numbers of \glspl{der} continues to fundamentally change the way that distribution networks operate. Combined siting/sizing and dispatch of these \glspl{der}, when considering the physical constraints of the underlying network, admits a non-convex, nonlinear \gls{minlp} problem form which scales very poorly in terms of solve time as the size of the network increases \cite{DeMel2024CompReform}. This paper extends previous solutions to this problem by De Mel et al. \cite{DeMel2024CompReform} through the incorporation of a distributed optimisation method, \gls{admm}, to solve these problems in a distributed manner, to decrease overall solution times and improve the scalability of the approach for large problems. A hybrid spatial/temporal decomposition process, exploiting the mathematical structure of the problem, is demonstrated for distribution networks with up to $55$ loads, with $120$ timepoints, corresponding to $1.9$ million variables and $4.5$ million constraints. Results are positive, showing computation speed-ups of up to $13$x under the assumption of solving all subproblems in parallel, yielding a maximum observed optimality gap of $0.61\%$. The results shown here demonstrate the utility of this approach, allowing larger network sizes to be solved in reasonable time frames. Like many \gls{admm} approaches however, empirical experiments to correctly set the algorithm parameters are still required to achieve optimal performance. Future work will focus on exploring efficient methods to determine these parameters, and how sensitive they are to changes in model parameters, for example the prices of different technologies. \\
The data used to produce the results in this paper is available at the Surrey Open Research repository \cite{Steven2026a} under a Creative Commons Attribution 4.0 International License.
    
    \section{Acknowledgements}
        Robert Steven is funded by the \gls{epsrc} through the Faculty of Engineering and Physical Sciences at the University of Surrey. For the purpose of open access, the authors have applied a Creative Commons attribution license (CC BY) to any Author Accepted Manuscript version arising from this submission.
    
    \printbibliography
    
    \begin{appendices}
\section{DES Model Parameters \& Variables}\label{app:parameters_variables}
\begin{landscape}
    \subsection{Parameters}
        \subsubsection{General}
            \begin{table}[H]
                \centering
                \caption{General Model Parameters}
                \begin{tabular}{|l|l|l|l|l|}
                    \hline
                    \textbf{Name} & \textbf{Value} & \textbf{Unit} & \textbf{Description} & \textbf{Source} \\
                    \hline
                    $n_{years}$ & $20$ & years & Number of years (for \gls{crf}) & Project \\
                    \hline
                    $n_{season}$ & $4$ & - & Number of seasons & - \\
                    \hline
                    $n_{day,s}$ & $\{92, 92, 91, 90\}$ & - & Number of days in each season & - \\
                    \hline
                    $\Delta t$ & $1$ & hour & Time interval & Project \\
                    \hline
                    $c^{water}$ & $1.0$ & $g/cm^3$ & Density of water & - \\
                    \hline
                    $\rho^{water}$ & $0.00116$ & $kWh/(kg\cdot ^{\circ}C)$ & Specific heat capacity of water & - \\
                    \hline
                    $EXPORT\_SEG\_TARIFF$ & $0.132$ & ${\pounds}/kWh$ & \Gls{seg} export price per kWh & Project \\
                    \hline
                    $R^{gas}$ & $0.02514$ & ${\pounds}/kWh$ & Natural gas price & \cite{DeMel2024CompReform} \\
                    \hline
                    Interest rate & 0.075 & - & Interest rate & \cite{DeMel2024CompReform} \\
                    \hline
                    $CRF$ & $0.0981$ & - & Capital Recovery Factor & \cite{DeMel2024CompReform} \\
                    \hline
                \end{tabular}
            \end{table}
    
        \subsubsection{Big M Values}
            \begin{table}[H]
                \centering
                \caption{Big M Values}
                \begin{tabular}{|l|l|l|l|}
                    \hline
                    \textbf{Name} & \textbf{Value} & \textbf{Description} \\
                    \hline
                    $M^{grid}$ & $100$ & Grid big M \\
                    \hline
                    $M^{batt,type}$ & $100$ & Battery type big M \\
                    \hline
                    $M^{batt,chg}$ & $100$ & Battery charge big M \\
                    \hline
                    $M^{boiler}$ & $100$ & Boiler big M \\
                    \hline
                    $M^{pump}$ & $100$ & Heat pump big M \\
                    \hline
                    $M^{tank}$ & $100$ & Hot water tank big M \\
                    \hline
                \end{tabular}
            \end{table}

        \subsubsection{Building Parameters}
            \begin{table}[H]
                \centering
                \caption{Building Parameters}
                \begin{tabular}{|l|l|l|l|l|}
                    \hline
                    \textbf{Name} & \textbf{Value} & \textbf{Unit} & \textbf{Description} & \textbf{Source} \\
                    \hline
                    $Vol^{building,avail}_h$ & $0.5$ & $m^3$ & Maximum building volume available for battery installation & \cite{DeMel2024CompReform} \\
                    \hline
                    $A^{roof\_max}_{h}$ & $35$ & $m^2$ & Area on building roof that can be used for \gls{pv} panels & \cite{DeMel2024CompReform} \\
                    \hline
                    $H^{building}_{h,t}$ & - & kWh & Building space heating load & Project (Section \ref{building_load_environ}) \\
                    \hline
                    $E^{building\_load}_{h,t}$ & - & kWh & Building electrical load & Project (Section \ref{building_load_environ}) \\
                    \hline
                    $PF$ & $0.95$ & $1$ & Building power factor & \cite{Schneider2018} \\
                    \hline
                    $T^{setpoint}$ & $20$ & $^{\circ}C$ & Temperature setpoint & \cite{DeMel2023HeatDecarb} \\
                    \hline
                \end{tabular}
            \end{table}
    
        \subsubsection{Environmental Parameters}
            \begin{table}[H]
                \centering
                \caption{Environmental Parameters}
                \begin{tabular}{|l|l|l|l|l|}
                    \hline
                    \textbf{Name} & \textbf{Value} & \textbf{Unit} & \textbf{Description} & \textbf{Source} \\
                    \hline
                    $Irradiance_{t}$ & - & $kW/m^2$ & Irradiance & \cite{DeMel2024CompReform} \\
                    \hline
                    $T^a_t$ & - & $^{\circ}C$ & Ambient temperature & \cite{Staffell2016,Pfenninger2016,RenewablesNinja} \\
                    \hline
                \end{tabular}
            \end{table}

        \subsubsection{Grid Interaction}
            \begin{table}[H]
                \centering
                \caption{Grid Parameters}
                \begin{tabular}{|l|l|l|l|l|}
                    \hline
                    \textbf{Name} & \textbf{Value} & \textbf{Unit} & \textbf{Description} & \textbf{Source} \\
                    \hline
                    $R^{grid,day}$ & $0.18$ & ${\pounds}/kWh$ & Grid electricity day tariff & \cite{DeMel2022b} \\
                    \hline
                    $R^{grid,night}$ & $0.08$ & ${\pounds}/kWh$ & Grid electricity night tariff & \cite{DeMel2022b} \\
                    \hline
                    $PG_{t}$ & - & ${\pounds}/kWh$ & Grid electricity tariff (night tariff between 00:00-07:00, day tariff otherwise) & \cite{DeMel2022b} \\
                    \hline
                    $S_{base}$ & $0.8$ & $MVA$ & Apparent power base & Project \\
                    \hline
                \end{tabular}
            \end{table}

        \subsubsection{PV}
            \begin{table}[H]
                \centering
                \caption{\Glsxtrshort{pv} Parameters}
                \begin{tabular}{|l|l|l|l|l|}
                    \hline
                    \textbf{Name} & \textbf{Value} & \textbf{Unit} & \textbf{Description} & \textbf{Source} \\
                    \hline			
                    $A^{PV\_panel}$ & $1.75$ & $m^2$ & Panel area & \cite{DeMel2024CompReform} \\
                    \hline
                    $\eta^{PV}$ & $0.18$ & $1$ & Panel efficiency & \cite{DeMel2024CompReform} \\
                    \hline
                    $Cap^{PV,panel}$ & $0.25$ & $kW$ & Panel capacity & \cite{DeMel2024CompReform} \\
                    \hline
                    $Cap^{PV,panel,max}$ & $5000$ & $kW$ & Maximum \gls{pv} panel capacity under \gls{seg} & \cite{SEG_Ofgem} \\
                    \hline
                    $R^{PV,capital}$ & $450.0$ & ${\pounds}/panel$ & Capital price per panel & \cite{DeMel2024CompReform} \\
                    \hline
                    $R^{PV,fixed\_op}$ & $12.5$ & ${\pounds}/kW-year$ & Fixed operating price & \cite{DeMel2024CompReform} \\
                    \hline
                \end{tabular}
            \end{table}
    
        \subsubsection{Batteries}
            Only a single lithium-ion battery type is considered here.
            \begin{table}[H]
                \centering
                \caption{Battery Parameters}
                \begin{tabular}{|l|l|l|l|l|}
                    \hline
                    \textbf{Name} & \textbf{Value} & \textbf{Unit} & \textbf{Description} & \textbf{Source} \\
                    \hline
                    $VED^{batt}_c$ & $148.37$ & $kWh/m^3$ & Volumetric energy density & \cite{FoxESSEP5Batt} \\
                    \hline
                    $SoC^{max}_{c}$ & $0.9$ & $1$ & Maximum \gls{soc} & \cite{DeMel2024CompReform} \\
                    \hline
                    $DoD^{max}_c$ & $0.9$ & $1$ & Maximum \gls{dod} & \cite{FoxESSEP5Batt} \\
                    \hline
                    $\eta^{batt,charge}_c$ & $0.97$ & $1$ & Charging efficiency & Calculated from round-trip efficiency of $0.95$ \cite{FoxESSEP5Batt} \\
                    \hline
                    $\eta^{batt,discharge}_c$ & $0.97$ & $1$ & Discharging efficiency & Calculated from round-trip efficiency of $0.95$ \cite{FoxESSEP5Batt} \\
                    \hline
                    $Chg^{batt,max\_rate}_c$ & $0.2$ & $1$ & Maximum charging rate & \cite{DeMel2024CompReform} \\
                    \hline
                    $Disch^{batt,max\_rate}_c$ & $0.2$ & $1$ & Maximum discharging rate & \cite{DeMel2024CompReform} \\
                    \hline
                    $R^{batt,capital}_c$ & $799$ & ${\pounds}/kWh$ & Capital price per kWh & \cite{Heatable} \\
                    \hline
                    $R^{batt,op}$ & $11$ & ${\pounds}/kWh-year$ & Fixed operating price & \cite{DeMel2024CompReform} \\
                    \hline
                \end{tabular}
            \end{table}

        \subsubsection{Boilers}
            \begin{table}[H]
                \centering
                \caption{Boiler Parameters}
                \begin{tabular}{|l|l|l|l|l|}
                    \hline
                    \textbf{Name} & \textbf{Value} & \textbf{Unit} & \textbf{Description} & \textbf{Source} \\
                    \hline
                    $\eta^{b}$ & $0.94$ & $1$ & Efficiency & \cite{DeMel2024CompReform} \\
                    \hline
                    $R^{boiler,capital}_b$ & $120$ & ${\pounds}/kW$ & Capital price & \cite{DeMel2024CompReform} \\
                    \hline
                \end{tabular}
            \end{table}

        \subsubsection{Heat Pumps}
            \begin{longtable}{|l|l|l|l|l|}
                \caption{Heat Pump Parameters} \\
                \hline
                \textbf{Name} & \textbf{Value} & \textbf{Unit} & \textbf{Description} & \textbf{Source} \\
                \hline
                \endfirsthead
                \multicolumn{5}{l}
                {\tablename\ \thetable\ -- \textit{Continued from previous page}} \\
                \hline
                \textbf{Name} & \textbf{Value} & \textbf{Unit} & \textbf{Description} & \textbf{Source} \\
                \hline
                \endhead
                \hline \multicolumn{5}{l}{\textit{Continued on next page}} \\
                \endfoot
                \hline
                \endlastfoot
                \hline
                $H^{p,max}_{t,hp}$ & - & kW & Maximum heating capacity & \cite{Mitsubishi_i} \\
                \hline
                $CoP^{pump}_{t,hp}$ & - & $1$ & \Gls{cop} & \cite{Mitsubishi_i} \\
                \hline
                $T^{pump,water\_supply}_{hp}$ & - & $^{\circ}C$ & Outlet water supply temperature & \cite{Mitsubishi_i} \\
                \hline
                $R^{pump,capital}_{hp}$ & - & £ & Capital price & \cite{CityPlumbing_a,CityPlumbing_b,CityPlumbing_c,CityPlumbing_d,CityPlumbing_e,SaturnSales_a} \\
                \hline
                $R^{pump,install}_{hp}$ & $3000$ & £ & Install price & Assumed \cite{Steven2025d} \\
                \hline
                $R^{pump,maint}_{hp}$ & $500$ & £ & Maintenance price & Assumed \\
            \end{longtable}
    
        \subsubsection{Hot Water Tanks}
            Two tanks from De Mel et al. \cite{DeMel2023HeatDecarb} are considered and two larger tanks are linearly extrapolated from these values.
            \begin{table}[H]
                \centering
                \caption{Hot Water Tank Parameters}
                \begin{tabular}{|p{1.8cm}|p{1.2cm}|p{1.2cm}|p{4cm}|p{4cm}|}
                    \hline
                    \textbf{Name} & \textbf{Value} & \textbf{Unit} & \textbf{Description} & \textbf{Source} \\
                    \hline
                    $T^{tank,min}$ & 49.0 & $^{\circ}C$ & Minimum tank water temperature & \cite{DeMel2023HeatDecarb} \\
                    \hline
                    $V^{tank}_w$ & - & $L$ & Volume & \cite{DeMel2023HeatDecarb} \\
                    \hline
                    $Loss^{tank}_w$ & - & kW & Heat losses & \cite{DeMel2023HeatDecarb} \\
                    \hline
                    $R^{tank,capital}_w$ & - & £ & Capital price & \cite{DeMel2023HeatDecarb} \\
                    \hline
                    $R^{tank,maint}_w$ & $0$ & £ & Maintenance price & Heat pump maintenance cost assumed to cover both heat pump and hot water tank \\
                    \hline
                \end{tabular}
            \end{table}
\end{landscape}

\begin{landscape}
    \subsection{Variables}
        \subsubsection{Type Key}
            \begin{table}[H]
                \centering
                \caption{Type Key}
                \begin{tabular}{|l|l|}
                    \hline
                    \textbf{Type} & \textbf{Description} \\
                    \hline
                    Bin & Binary variable, can only be $0$ or $1$ \\
                    \hline
                    NNR & Non-negative real numbers ($0$ and positive numbers) \\
                    \hline
                    R & All real numbers (negative, $0$ and positive numbers) \\
                    \hline
                \end{tabular}
            \end{table}

        \subsubsection{Sets}
            \begin{table}[H]
                \centering
                \caption{Sets}
                \begin{tabular}{|l|l|}
                    \hline
                    \textbf{Type} & \textbf{Description} \\
                    \hline
                    $\mathcal{H}$ & Loads \\
                    \hline
                    $\mathcal{T}$ & Timepoints \\
                    \hline
                    $\mathcal{S}$ & Seasons \\
                    \hline
                    $\mathcal{SR}$ & Robust seasons \\
                    \hline
                    $\mathcal{C}$ & Batteries \\
                    \hline
                    $\mathcal{B}$ & Boilers \\
                    \hline
                    $\mathcal{HP}$ & Heat pumps \\
                    \hline
                    $\mathcal{W}$ & Hot water tanks \\
                    \hline
                \end{tabular}
            \end{table}

        \subsubsection{Energy Balance}
            \begin{longtable}{|l|l|l|l|}
                \caption{Energy Balance Variables} \\
                \hline
                \textbf{Name} & \textbf{Type} & \textbf{Unit} & \textbf{Description} \\
                \hline
                \endfirsthead
                \multicolumn{4}{l}
                {\tablename\ \thetable\ -- \textit{Continued from previous page}} \\
                \hline
                \textbf{Name} & \textbf{Type} & \textbf{Unit} & \textbf{Description} \\
                \hline
                \endhead
                \hline \multicolumn{4}{l}{\textit{Continued on next page}} \\
                \endfoot
                \hline
                \endlastfoot
                $E^{load}_{h,t}$ & NNR & kWh & Total electrical load for each building \\
            \end{longtable}
    
        \subsubsection{Model Costs \& Income}
            \begin{table}[H]
                \centering
                \caption{Model Costs \& Income Variables}
                \begin{tabular}{|l|l|l|l|}
                    \hline
                    \textbf{Name} & \textbf{Type} & \textbf{Unit} & \textbf{Description} \\
                    \hline
                    $C^{income,export}_{s}$ & NNR & £ & Income from exporting electricity to National Grid under \gls{seg}\\
                    \hline
                \end{tabular}
            \end{table}
    
        \subsubsection{Grid Interaction}
            \begin{table}[H]
                \centering
                \caption{Grid Interaction Variables}
                \begin{tabular}{|l|l|l|l|}
                    \hline
                    \textbf{Name} & \textbf{Type} & \textbf{Unit} & \textbf{Description} \\
                    \hline
                    $X_{t}$ & Bin & - & 1 if electricity sold to the grid \\
                    \hline
                    $E^{grid}_{h,t}$ & NNR & kWh & Electricity imported from the wider grid \\
                    \hline
                    $E^{grid,load}_{h,t}$ & NNR & kWh & Electricity imported from the wider grid for building load \\
                    \hline
                    $E^{grid,charge}_{h,t,c}$ & NNR & kWh & Electricity imported from the wider grid for battery charging \\
                    \hline
                    $C^{OPEX,grid}_{s}$ & NNR & £ & Grid electricity cost \\
                    \hline
                    $P^{inject}_{h,\psi,t}$ & R & kW & Real power injected into the network \\
                    \hline
                    $Q^{inject}_{h,\psi,t}$ & R & kVAr & Reactive power injected into the network \\
                    \hline
                \end{tabular}
            \end{table}
    
        \subsubsection{PV}
            \begin{table}[H]
                \centering
                \caption{\Glsxtrshort{pv} Variables}
                \begin{tabular}{|l|l|l|l|}
                    \hline
                    \textbf{Name} & \textbf{Type} & \textbf{Unit} & \textbf{Description} \\
                    \hline
                    $PV^{panels}_{h}$ & NNR & - & Number of installed \gls{pv} panels \\
                    \hline
                    $E^{PV,used}_{h,t}$ & NNR & kWh & Energy generated by \glspl{pv} used to meet building load \\
                    \hline
                    $E^{PV,sold}_{h,t}$ & NNR & kWh & Energy generated by \glspl{pv} sold back to grid \\
                    \hline
                    $C^{PV,CAPEX}_{s}$ & NNR & £ & \Gls{pv} \gls{capex} \\
                    \hline
                    $C^{PV,OPEX}_{s}$ & NNR & £ & \Gls{pv} \gls{opex} \\
                    \hline
                \end{tabular}
            \end{table}
    
        \subsubsection{Batteries}
            \begin{table}[H]
                \centering
                \caption{Battery Variables}
                \begin{tabular}{|l|l|l|l|}
                    \hline
                    \textbf{Name} & \textbf{Type} & \textbf{Unit} & \textbf{Description} \\
                    \hline
                    $W_{h,c}$ & Bin & - & 1 if battery selected \\
                    \hline
                    $Q^{batt}_{h,t,c}$ & Bin & - & 1 if battery charging \\
                    \hline
                    $Cap^{batt}_{h,c}$ & NNR & kWh & Installed battery capacity \\
                    \hline
                    $Vol^{batt}_{h,c}$ & NNR & $m^3$ & Installed battery volume \\
                    \hline
                    $E^{batt,stored}_{h,t,c}$ & NNR & kWh & Energy stored \\
                    \hline
                    $E^{batt,charge}_{h,t,c}$ & NNR & kWh & Total energy charged \\
                    \hline
                    $E^{batt,discharge}_{h,t,c}$ & NNR & kWh & Total energy discharged \\
                    \hline
                    $E^{PV,charge}_{h,t,c}$ & NNR & kWh & Energy from \gls{pv} used to charge battery \\
                    \hline
                    $C^{batt,CAPEX}_{s}$ & NNR & £ & Battery capital cost \\
                    \hline
                    $C^{batt,OPEX}_{s}$ & NNR & £ & Battery operating cost \\
                    \hline
                \end{tabular}
            \end{table}

        \subsubsection{Boilers}
            \begin{table}[H]
                \centering
                \caption{Boiler Variables}
                \begin{tabular}{|l|l|l|l|}
                    \hline
                    \textbf{Name} & \textbf{Type} & \textbf{Unit} & \textbf{Description} \\
                    \hline
                    $B_{h}$ & Bin & - & 1 if boiler selected \\
                    \hline
                    $H^{b}_{h,t}$ & NNR & kW & Heat energy output \\
                    \hline
                    $H^{b,max}_{h}$ & NNR & kW & Maximum boiler heat energy output \\
                    \hline
                    $C^{boiler,CAPEX}_{s}$ & NNR & £ & Boiler capital cost \\
                    \hline
                    $C^{boiler,OPEX}_{s}$ & NNR & £ & Boiler operating cost \\
                    \hline
                \end{tabular}
            \end{table}

        \subsubsection{Heat Pumps}
            \begin{table}[H]
                \centering
                \caption{Heat Pump Variables}
                \begin{tabular}{|l|l|l|l|}
                    \hline
                    \textbf{Name} & \textbf{Type} & \textbf{Unit} & \textbf{Description} \\
                    \hline
                    $P_{h,hp}$ & Bin & - & 1 if heat pump selected \\
                    \hline
                    $H^p_{h,t,hp,w}$ & NNR & kWh & Heat energy output \\
                    \hline
                    $E^{H\_p}_{h,t,hp}$ & NNR & kWh & Electrical energy input \\
                    \hline
                    $C^{pump,CAPEX}_{s}$ & NNR & £ & Pump capital cost \\
                    \hline
                    $C^{pump,OPEX}_{s}$ & NNR & £ & Pump operating cost \\
                    \hline
                \end{tabular}
            \end{table}
    
        \subsubsection{Hot Water Tanks}
            \begin{table}[H]
                \centering
                \caption{Hot Water Tank Variables}
                \begin{tabular}{|l|l|l|l|}
                    \hline
                    \textbf{Name} & \textbf{Type} & \textbf{Unit} & \textbf{Description} \\
                    \hline
                    $Z_{h,w}$ & Bin & - & 1 if hot water tank selected \\
                    \hline
                    $T^{tank}_{h,t,w}$ & NNR & $^{\circ}C$ & Internal temperature of tank \\
                    \hline
                    $H^{tank\_internal}_{h,t,w}$ & NNR & kWh & Internal heat of tank \\
                    \hline
                    $H^{tank,charge}_{h,t,w}$ & NNR & kWh & Heat charged to tank \\
                    \hline
                    $H^{tank,discharge}_{h,t,w}$ & NNR & kWh & Heat discharged (delivered) by tank \\
                    \hline
                    $C^{tank,CAPEX}_{s}$ & NNR & £ & Tank capital cost \\
                    \hline
                    $C^{tank,OPEX}_{s}$ & NNR & £ & Tank operating cost \\
                    \hline
                \end{tabular}
            \end{table}
\end{landscape}

\section{DES Modelling Component Equations}\label{app:modelling_components}
    \subsection{Objective Function}
\subsubsection{Economic Cost Function}
    The total economic cost to minimise.
    \begin{equation}
        TAC=DES\_SCALE\_FACTOR\left( C^{CAPEX} + C^{OPEX} - C^{income} \right)
    \end{equation}

\subsubsection{Capital Cost}
    \Gls{capex} costs.
    \begin{equation}
        C^{CAPEX} = \sum_{s\in\mathcal{S}} \left[ C^{CAPEX,PV}_s + C^{CAPEX,batt}_s + C^{CAPEX,boiler}_s + C^{CAPEX,pump}_s + C^{CAPEX,tank}_s \right]
    \end{equation}

\subsubsection{Operational Cost}
    \Gls{opex} costs. Note that for brevity, seasonal indexes on variables other than those making up each cost component are omitted.
    \begin{equation}
        C^{OPEX} = \sum_{s\in\mathcal{S}} \left[ C^{OPEX,PV}_s + C^{OPEX,batt}_s + C^{OPEX,boiler}_s + C^{OPEX,pump}_s + C^{OPEX,tank}_s + C^{OPEX,grid}_s \right]
    \end{equation}

\subsubsection{Income}
    Income from \gls{seg}.
    \begin{equation}
        C^{income} = \sum_{s\in\mathcal{S}} \left[ C^{income,export}_{s} \right]
    \end{equation}

    \subsection{Energy Balances}
\subsubsection{Heating Balance}
    Heat energy balance (load and generation) for each building.
    \begin{equation}
        H^{building}_{h,t} = \left(H^{boiler}_{h,t} \cdot \Delta t\right) + \sum_{w\in\mathcal{W}} \left[ H^{tank\_disch}_{h,t,w} \right], \forall h\in\mathcal{H}, t\in\mathcal{T}
    \end{equation}
	
\subsubsection{Electrical Load}
    Electrical load for each building.
    \begin{equation}
        E^{load}_{h,t} = E^{building\_load}_{h,t} + \sum_{hp\in\mathcal{HP}} \left[ E^{H\_p}_{h,t,hp} \right], \forall h\in\mathcal{H}, t\in\mathcal{T}
    \end{equation}

\subsubsection{Electrical Balance}
    Electrical energy balance (load and generation) for each building.
    \begin{equation}
        E^{load}_{h,t} = E^{grid,load}_{h,t} + E^{PV,used}_{h,t} + \sum_{c\in\mathcal{C}} \left[ E^{batt,discharge}_{h,t,c} \right], \forall h\in\mathcal{H}, t\in\mathcal{T}
    \end{equation}
    
    \subsection{Photovoltaics}
Modelling equations from \cite{DeMel2024CompReform}.

\subsubsection{PV Generation}
    The number of panels, multiplied by their area, the irradiance on each panel and the efficiency of the panel gives the generated \gls{pv} power that can be used to satisfy building load, sold to the grid or used to charge local battery.
    \begin{equation} 
        E^{PV,used}_{h,t} + E^{PV,sold}_{h,t} + \sum_{c\in\mathcal{C}} \left[ E^{PV,charge}_{h,t,c} \right] \leq PV^{panels}_{h} \cdot A^{PV\_panel} \cdot Irradiance_{t} \cdot \eta^{PV} \cdot \Delta t, \forall h\in\mathcal{H}, t\in\mathcal{T}
    \end{equation}
	
\subsubsection{PV Maximum Generation Capacity}
    Generated \gls{pv} power is less than or equal to the number of installed panels multiplied by the maximum capacity of the panel.
    \begin{equation}
        E^{PV,used}_{h,t} + E^{PV,sold}_{h,t} + \sum_{c\in\mathcal{C}} \left[ E^{PV,charge}_{h,t,c} \right] \leq PV^{panels}_{h} \cdot Cap^{PV,panel} \cdot \Delta t, \forall h\in\mathcal{H}, t\in\mathcal{T}
    \end{equation}

\subsubsection{Maximum Roof Area}
    The total area of the installed panels is less than or equal to the area of available rooftop space on each building.
    \begin{equation}
        PV^{panels}_{h} \cdot A^{PV\_panel} \leq A^{roof\_max}_{h}, \forall h\in\mathcal{H}
    \end{equation}

\subsubsection{PV Capacity Limitation}
    Total installed \gls{pv} capacity limited by \gls{seg} regulations.
    \begin{equation}
        PV^{panels}_{h} \cdot Cap^{PV,panel} \leq Cap^{PV,panel,max}, \forall h\in\mathcal{H}
    \end{equation}

\subsubsection{PV CAPEX}
    Capital expenditure for all installed \gls{pv}, across all buildings.
    \begin{equation} 
        C^{PV,CAPEX}_{s} = \sum_{h\in\mathcal{H}} \left[ R^{PV,capital} \cdot PV^{panels}_{h} \cdot CRF \right] \cdot \frac{1}{n_{season}}, \forall s\in\mathcal{S}
    \end{equation}

\subsubsection{PV OPEX}
    Operational expenditure for all installed \gls{pv}, across all buildings.
    \begin{equation}
        C^{PV,OPEX}_{s} = \sum_{h\in\mathcal{H}} \left[ PV^{panels}_h \cdot  R^{PV,fixed\_op} \cdot \frac{1}{365} \cdot n_{day,s} \cdot Cap^{PV\_panel} \right], \forall s\in\mathcal{S}
    \end{equation}

\subsubsection{PV Panel Seasonal Linking Constraint}
    \Gls{pv} installation decisions are the same across all seasons.
    \begin{equation}
        \begin{cases}
            PV^{panels}_{s,h} & \text{if } s=1 \\
            PV^{panels}_{s,h} = PV^{panels}_{s-1,h} & \text{otherwise}
        \end{cases}
    \end{equation}

\subsubsection{PV Panel Robust Linking Constraint}
    \Gls{pv} installation decisions are the same for robust seasons as they are for normal seasons.
    \begin{equation}
        PV^{panels}_{s_r,h} = PV^{panels}_{s=1,h}, \forall s_r\in\mathcal{SR}
    \end{equation}

    \subsection{Battery Storage}
Modelling equations from \cite{DeMel2024CompReform}.

\subsubsection{Battery Type}
    Only install a maximum of 1 battery option at each building.
    \begin{equation} 
         \sum_{c\in\mathcal{C}} \left[ W_{h,c} \right] \leq 1, \forall h\in\mathcal{H}, c\in\mathcal{C}
    \end{equation}

\subsubsection{Battery Type BigM}
    Can only define a non-zero battery capacity if the battery is installed.
    \begin{equation}
         Cap^{batt}_{h,c} \leq M^{batt,type} \cdot W_{h,c}, \forall h\in\mathcal{H}, c\in\mathcal{C}
    \end{equation}

\subsubsection{Installed Battery Capacity}
    The installed battery capacity is equal to the installed battery volume multiplied by the battery \gls{ved}.
    \begin{equation}
        Cap^{batt}_{h,c} = Vol^{batt}_{h,c} \cdot VED^{batt}_c, \forall h\in\mathcal{H}, c\in\mathcal{C}
    \end{equation}

\subsubsection{Battery Volume Limit}
    The total installed battery volume is less than or equal to the available installed volume in each building.
    \begin{equation}
        \sum_{c\in\mathcal{C}} \left[ Vol^{batt}_{h,c} \right] \leq Vol^{building,avail}_h, \forall h\in\mathcal{H}
    \end{equation}

\subsubsection{Battery Capacity 1}
    The maximum \gls{soc} sets the upper limit on the energy stored in the battery.
    \begin{equation}
         E^{batt,stored}_{h,t,c} \leq Cap^{batt}_{h,c} \cdot SoC^{max}_{c}, \forall h\in\mathcal{H}, t\in\mathcal{T}, c\in\mathcal{C}
    \end{equation}

\subsubsection{Battery Capacity 2}
    The maximum \gls{dod} sets the lower limit on the energy stored in the battery.
    \begin{equation}
        E^{batt,stored}_{h,t,c} \geq Cap^{batt}_{h,c} \cdot (1 - DoD^{max}_c), \forall h\in\mathcal{H}, t\in\mathcal{T}, c\in\mathcal{C}
    \end{equation}

\subsubsection{Battery Storage Balance}
    Equation defining the energy stored in the battery at each timepoint, from energy charged into the battery.
    \begin{multline}
        E^{batt,stored}_{h,t,c} =
        \begin{cases}
            E^{batt,stored}_{h,t-1,c} + \left( E^{batt,charge}_{h,t,c} \cdot \eta^{batt,charge}_c \right) - \\ \left( \frac{E^{batt,discharge}_{h,t,c}}{\eta^{batt,discharge}_c} \right) & \text{if } t>1 \\
            \left( E^{batt,charge}_{h,t,c} \cdot \eta^{batt,charge}_c \right) - \left( \frac{E^{batt,discharge}_{h,t,c}}{\eta^{batt,discharge}_c} \right) & \text{otherwise}
        \end{cases}, \\ \forall h\in\mathcal{H}, t\in\mathcal{T}, c\in\mathcal{C}
    \end{multline}

\subsubsection{Battery Discharge Condition}
    Equation defining the energy stored in the battery at each timepoint, from energy discharged from the battery.
    \begin{equation}
        \frac{E^{batt,discharge}_{h,t,c}}{\eta^{batt,discharge}_c} \leq E^{batt,stored}_{h,t-1,c}\quad \text{if } t>1, \forall h\in\mathcal{H}, t\in\mathcal{T}, c\in\mathcal{C}
    \end{equation}

\subsubsection{Battery Start and End SoC}
    The energy in the battery is kept the same at the start and end of the representational seasonal time period.
    \begin{equation}
         E^{batt,stored}_{h,t=1,c} = E^{batt,stored}_{h,t=t_{end},c}, \forall h\in\mathcal{H}, c\in\mathcal{C}
    \end{equation}

\subsubsection{Battery Charge Limitation}
    Limit on the maximum charging rate of the battery.
    \begin{equation} 
         E^{batt,charge}_{h,t,c} \cdot \eta^{batt,charge}_c \leq Cap^{batt}_{h,c} \cdot Chg^{batt,max\_rate}_c, \forall h\in\mathcal{H}, t\in\mathcal{T}, c\in\mathcal{C}
    \end{equation}

\subsubsection{Battery Discharge Limitation}
    Limit on the maximum discharge rate of the battery.
    \begin{equation}
         \frac{E^{batt,discharge}_{h,t,c}}{\eta^{batt,discharge}_c} \leq Cap^{batt}_{h,c} \cdot Disch^{batt,max\_rate}_c, \forall h\in\mathcal{H}, t\in\mathcal{T}, c\in\mathcal{C}
    \end{equation}

\subsubsection{Battery Charge}
    Charge for installed batteries can come from local \gls{pv} or the grid.
    \begin{equation}
        E^{batt,charge}_{h,t,c} = E^{PV,charge}_{h,t,c} + E^{grid,charge}_{h,t,c}, \forall h\in\mathcal{H}, t\in\mathcal{T}, c\in\mathcal{C}
    \end{equation}

\subsubsection{Battery Charge BigM}
    If the battery is charging in a given time period then it cannot also discharge.
    \begin{equation}
         E^{batt,charge}_{h,t,c} \leq M^{batt,chg} \cdot Q^{batt}_{h,t,c}, \forall h\in\mathcal{H}, t\in\mathcal{T}, c\in\mathcal{C}
    \end{equation}

\subsubsection{Battery Discharge BigM}
    If a battery is discharging in a given time period then it cannot also charge.
    \begin{equation}
         E^{batt,discharge}_{h,t,c} \leq M^{batt,chg} \cdot \left(1 - Q^{batt}_{h,t,c} \right), \forall h\in\mathcal{H}, t\in\mathcal{T}, c\in\mathcal{C}
    \end{equation}

\subsubsection{Battery CAPEX}
    Capital expenditure for all installed batteries, across all buildings.
    \begin{equation}
        C^{batt,CAPEX}_{s} = \sum_{h\in\mathcal{H}} \sum_{c\in\mathcal{C}} \left[ Cap^{batt}_{h,c} \cdot R^{batt,capital}_c \cdot CRF \right] \cdot \frac{1}{n_{season}}, \forall s\in\mathcal{S}
    \end{equation}

\subsubsection{Battery OPEX}
    Operational expenditure for all installed batteries, across all buildings.
    \begin{equation}
        C^{batt,OPEX}_{s} = \sum_{h\in\mathcal{H}} \sum_{c\in\mathcal{C}} \left[ Cap^{batt}_{h,c} \cdot R^{batt,op}_c \right] \cdot \frac{1}{365} \cdot n_{day,s}, \forall s\in\mathcal{S}
    \end{equation}

\subsubsection{Battery Seasonal Linking Constraint}
    Battery installation decisions are the same across all seasons.
    \begin{equation}
        \begin{cases}
            Cap^{batt}_{s,h,c} & \text{if } s=1 \\
            Cap^{batt}_{s,h,c} = Cap^{batt}_{s-1,h,c} & \text{otherwise} 
        \end{cases}
    \end{equation}

\subsubsection{Battery Robust Linking Constraint}
    Battery installation decisions are the same for robust seasons as they are for normal seasons.
    \begin{equation}
        Cap^{batt}_{s_r,h,c} = Cap^{batt}_{s=1,h,c}, \forall s_r\in\mathcal{SR}
    \end{equation}
        
    \subsection{Grid Interaction}
Modelling equations from \cite{DeMel2024CompReform}.

\subsubsection{Grid Electricity Usage}
    \begin{equation}
        E^{grid}_{h,t} = E^{grid,load}_{h,t} + \sum_{c\in\mathcal{C}} \left[ E^{grid,charge}_{h,t,c} \right], \forall h\in\mathcal{H}, t\in\mathcal{T}
    \end{equation}

\subsubsection{Grid Electricity Usage for Local Load}
    Electricity purchased from the grid to meet local building load must be less than or equal to that local building load.
    \begin{equation}
        E^{grid,load}_{h,t} \leq E^{load}_{h,t}, \forall h\in\mathcal{H}, t\in\mathcal{T}
    \end{equation}

\subsubsection{Grid Electricity Purchase}
    \begin{equation}
        E^{grid}_{h,t} \leq M^{grid} \cdot \left(1 - X_{h,t}\right), \forall h\in\mathcal{H}, t\in\mathcal{T}
    \end{equation}

\subsubsection{Grid Electricity Sale}
    Electricity sold to the grid (by each building) comes from generated \gls{pv} electricity.
    \begin{equation}
        E^{PV,sold}_{h,t} \leq M^{grid} \cdot X_{h,t}, \forall h\in\mathcal{H}, t\in\mathcal{T}
    \end{equation}
	
\subsubsection{Cost of Purchasing Electricity}
    Total cost of purchasing electricity from the grid.
    \begin{equation}
        C^{OPEX,grid}_{s} = \sum_{h\in\mathcal{H}} \sum_{t\in\mathcal{T}} \left[ E^{grid}_{h,t} \cdot PG_{t} \right] \cdot n_{day,s}, \forall s\in\mathcal{S}
    \end{equation}

\subsubsection{Real Power Injection}
    Real power injection into the network.
    \begin{equation}
        P^{inject}_{h,\psi_{load},t} = 1*10^{-3} \cdot \left(E^{PV,sold}_{h,t} - E^{grid}_{h,t} \right) \cdot \frac{1}{S_{base}} \cdot \frac{1}{\Delta t}, \forall h\in\mathcal{H}, t\in\mathcal{T}
    \end{equation}

\subsubsection{Reactive Power Injection}
    Reactive power injection into the network.
    \begin{equation}
        Q^{inject}_{h,\psi_{load},t} = -1*10^{-3} \cdot \sqrt{\left(E^{building\_load}_{h,t}\right)^2\cdot\left(\frac{1}{PF^{2}}-1\right)} \cdot \frac{1}{S_{base}} \cdot \frac{1}{\Delta t}, \forall h\in\mathcal{H}, t\in\mathcal{T}
    \end{equation}

    \subsection{Income}
Modelling equations from \cite{DeMel2024CompReform}.

\subsubsection{Export Income}
    Income from exporting power to grid under the \gls{seg} rules.
    \begin{equation}
        C^{income,export}_{s} = \sum_{h\in\mathcal{H}} \sum_{t\in\mathcal{T}} \left[ E^{PV,sold}_{h,t} \cdot EXPORT\_SEG\_TARIFF \right] \cdot n_{day,s}, \forall s\in\mathcal{S}
    \end{equation}

    \subsection{Boiler}
Modelling equations from \cite{DeMel2024CompReform}.

\subsubsection{Boiler BigM}
    Boiler can only output heat if it is installed.
    \begin{equation}
        H^{b,max}_{h} \leq M^{boiler} \cdot B_{h}, \forall h\in\mathcal{H}, t\in\mathcal{T}, b\in\mathcal{B}
    \end{equation}

\subsubsection{Maximum Boiler Capacity}
    Boiler maximum heating capacity.
    \begin{equation}
        H^{b}_{h,t} \leq H^{b,max}_{h}, \forall h\in\mathcal{H}, t\in\mathcal{T}, b\in\mathcal{B}
    \end{equation}

\subsubsection{Boiler CAPEX}
    \Gls{capex} for all installed boilers across all buildings.	
    \begin{equation}
        C^{boiler,CAPEX}_{s} = \sum_{h\in\mathcal{H}} \sum_{b\in\mathcal{B}} \left[ R^{boiler,capital} \cdot H^{b,max}_{h} \cdot CRF \right] \cdot \frac{1}{n_{season}}, \forall s\in\mathcal{S}
    \end{equation}

\subsubsection{Boiler OPEX}
    \Gls{opex} for all installed boilers, across all buildings.
    \begin{equation}
        C^{boiler,OPEX}_{s} = \sum_{h\in\mathcal{H}} \sum_{b\in\mathcal{B}} \sum_{t\in\mathcal{T}} \left[ H^{b}_{h,t} \cdot \Delta t \cdot \frac{R^{gas}}{\eta^{b}} \right] \cdot n_{day,s}, \forall s\in\mathcal{S}
    \end{equation}

\subsubsection{Boiler Seasonal Linking Constraint}
    Boiler installation decisions are the same across all seasons.
    \begin{equation}
        \begin{cases}
            H^{s,b,max}_{h} & \text{if } s=1 \\
            H^{s,b,max}_{h} = H^{s-1,b,max}_{h} & \text{otherwise}
        \end{cases}, \forall h\in\mathcal{H}
    \end{equation}

\subsubsection{Boiler Robust Linking Constraint}
    Boiler installation decisions are the same for robust seasons as they are for normal seasons.
    \begin{equation}
            H^{s_r,b,max}_{h} = H^{s=1,b,max}_{h}, \forall s_r\in\mathcal{SR}, h\in\mathcal{H}
    \end{equation}
    
    \subsection{Heat Pump}
Modelling equations from \cite{DeMel2023HeatDecarb}.

\subsubsection{Heat Pump Limitation}
    Only install a maximum of 1 heat pump options at each building.
    \begin{equation} 
         \sum_{hp\in\mathcal{HP}} \left[ P_{h,hp} \right] \leq 1,\forall h\in\mathcal{H}
    \end{equation}

\subsubsection{Heat Pump BigM}
    Heat pump can only output heat if it is installed.
    \begin{equation} 
         H^p_{h,t,hp,w} \leq M^{pump} \cdot P_{h,hp},\forall h\in\mathcal{H}, t\in\mathcal{T}, hp\in\mathcal{HP}, w\in\mathcal{W}
    \end{equation}

\subsubsection{Heat Pump Electrical Load}
    Heat pump electrical load is equal to its heat output divided by its \gls{cop}. E.g. for a heat pump outputting 3 kW of heat, with a \gls{cop} of 3 it will only require 1 kW of electrical power.
    \begin{equation} 
         E^{H\_p}_{h,t,hp} = \frac{\sum_{w\in\mathcal{W}} \left[ H^p_{h,t,hp,w} \right] }{CoP^{pump}_{t,hp}}, \forall h\in\mathcal{H}, t\in\mathcal{T}, hp\in\mathcal{HP}
    \end{equation}

\subsubsection{Heat Pump Maximum Capacity}
    The heat output of the heat pump is equal to or less than its maximum heating capacity.
    \begin{equation} 
         \sum_{w\in\mathcal{W}} \left[ H^p_{h,t,hp,w} \right] \leq H^{p,max}_{t,hp} \cdot \Delta t,\forall h\in\mathcal{H}, t\in\mathcal{T}, hp\in\mathcal{HP}
    \end{equation}

\subsubsection{Heat Pump CAPEX}
    Capital expenditure for all installed heat pumps, across all buildings.
    \begin{equation}
        C^{pump,CAPEX}_{s} = \sum_{h\in\mathcal{H}} \sum_{hp\in\mathcal{HP}} \Bigg[ P_{h,hp} \cdot \left( R^{pump,capital}_{hp} + R^{pump,install}_{hp} \right) \cdot CRF \Bigg] \cdot \frac{1}{n_{season}}, \forall s\in\mathcal{S}
    \end{equation}

\subsubsection{Heat Pump OPEX}
    Operational expenditure for all installed heat pumps, across all buildings.
    \begin{equation}
        C^{pump,OPEX}_{s} = \sum_{h\in\mathcal{H}} \sum_{hp\in\mathcal{HP}} \left[ P_{h,hp} \cdot R^{pump,maint}_{hp} \right] \cdot \frac{1}{n_{season}}, \forall s\in\mathcal{S}
    \end{equation}

\subsubsection{Heat Pump Seasonal Linking Constraint}
    Heat pump installation decisions are the same across all seasons.
    \begin{equation}
        \begin{cases}
            P_{s,h,hp} & \text{if } s=1 \\
            P_{s,h,hp} = P_{s-1,h,hp} & \text{otherwise}
        \end{cases}, \forall h\in\mathcal{H}, hp\in\mathcal{HP}
    \end{equation}

\subsubsection{Heat Pump Robust Linking Constraint}
    Heat pump installation decisions are the same for robust seasons as they are for normal seasons.
    \begin{equation}
            P_{s_r,h,hp} = P_{s=1,h,hp}, \forall s_r\in\mathcal{SR}, h\in\mathcal{H}, hp\in\mathcal{HP}
    \end{equation}
        
    \subsection{Heating Technologies}
\subsubsection{Simultaneous Heating Technologies}
    Can only install either a boiler or a heat pump at each load.
    \begin{equation}
        P_{h,hp} + B_{h,b} \leq 1, \forall h\in\mathcal{H}
    \end{equation}
    
    \subsection{Hot Water Tank}
Modelling equations from \cite{DeMel2023HeatDecarb}.

\subsubsection{Hot Water Tank Limitation}
    Only install a maximum of 1 hot water tank options at each building.
    \begin{equation} 
         \sum_{w\in\mathcal{W}} \left[ Z_{h,w} \right] \leq 1, \forall h\in\mathcal{H}
    \end{equation} 

\subsubsection{Hot Water Tank Heat Charge}
    Heat energy charged into the hot water tank.
    \begin{equation} 
         H^{tank,charge}_{h,t,w} = \sum_{hp\in\mathcal{HP}} \left[ H^p_{h,t,hp,w} \right], \forall h\in\mathcal{H}, t\in\mathcal{T}, w\in\mathcal{W}
    \end{equation} 

\subsubsection{Hot Water Tank BigM Limitation}
    Heat can only be charged into the hot water tank if it is installed.
    \begin{equation} 
         H^{tank,charge}_{h,t,w} \leq M^{tank} \cdot Z_{h,w}, \forall h\in\mathcal{H}, t\in\mathcal{T}, w\in\mathcal{W}
    \end{equation} 

\subsubsection{Hot Water Tank Minimum Temperature}
    Temperature of the water in the hot water tank must stay above the minimum tank temperature.
    \begin{equation} 
         T^{tank,min} \cdot Z_{h,w} \leq T^{tank}_{h,t,w}, \forall h\in\mathcal{H}, t\in\mathcal{T}, w\in\mathcal{W}
    \end{equation} 

\subsubsection{Hot Water Tank Temperature}
    Temperature of the water in the hot water tank.
    \begin{equation}
        T^{tank}_{h,t,w} = \frac{H^{tank\_internal}_{h,t,w}}{c^{water} \cdot \rho^{water} \cdot V^{tank}_w} + T^{setpoint} \cdot Z_{h,w}, \forall h\in\mathcal{H}, t\in\mathcal{T}, w\in\mathcal{W}
    \end{equation}

\subsubsection{Hot Water Tank Maximum Temperature}
    Maximum temperature of the water that can be charged into the hot water tank is set by the outlet temperature of the heat pump that is supplying the hot water.
    \begin{equation}
        T^{tank}_{h,t,w} \leq \sum_{hp\in\mathcal{HP}} \left[ T^{pump,water\_supply}_{hp}\cdot P_{h,hp} \right],\forall h\in\mathcal{H}, t\in\mathcal{T}, w\in\mathcal{W}
    \end{equation}

\subsubsection{Hot Water Tank Heat}
    Heat energy of the water in the hot water tank.
    \begin{multline}
        H^{tank\_internal}_{h,t,w} = 
        \begin{cases}
            H^{tank\_internal}_{h,t-1,w} + \bigg( H^{tank,charge}_{h,t,w} - \Big( H^{tank,discharge}_{h,t,w} + \\ Z_{h,w} \cdot Loss^{tank}_w \cdot \Delta t \Big) \bigg) & \text{if } t>1 \\
            Z_{h,w} \cdot \left( T^{tank,min} - T^{setpoint} \right) \cdot \left(c^{water} \cdot \rho^{water} \cdot V^{tank}_w \right) + \\ \left( H^{tank,charge}_{h,t,w} - \left( H^{tank,discharge}_{h,t,w} + Z_{h,w} \cdot Loss^{tank}_w \cdot \Delta t \right) \right) & \text{otherwise}
        \end{cases}, \\ \forall h\in\mathcal{H}, t\in\mathcal{T}, w\in\mathcal{W}
    \end{multline}

\subsubsection{Hot Water Tank Heat Timebounds}
    The heat in the hot water tank is kept the same at the start and end of the representational seasonal time period.
    \begin{equation} 
         H^{tank\_internal}_{h,t=1,w} = H^{tank\_internal}_{h,t=t_{end},w}, \forall h\in\mathcal{H}, t\in\mathcal{T}, w\in\mathcal{W}
    \end{equation} 

\subsubsection{Hot Water Tank CAPEX}
    \Gls{capex} for all installed hot water tanks, across all buildings.
    \begin{equation}
        C^{tank,CAPEX}_{s} = \sum_{h\in\mathcal{H}} \sum_{w\in\mathcal{W}} \Bigg[ Z_{h,w} \cdot R^{tank,capital}_w \cdot CRF \Bigg] \cdot \frac{1}{n_{season}}, \forall s\in\mathcal{S}
    \end{equation}

\subsubsection{Hot Water Tank OPEX}
    \Gls{opex} for all installed hot water tanks, across all buildings.
    \begin{equation}
        C^{tank,OPEX}_{s} = \sum_{h\in\mathcal{H}} \sum_{w\in\mathcal{W}} \left[ Z_{h,w} \cdot R^{tank,maint}_w \right] \cdot \frac{1}{n_{season}}, \forall s\in\mathcal{S}
    \end{equation}

\subsubsection{Hot Water Tank Seasonal Linking Constraint}
    Hot water tank installation decisions are the same across all seasons.
    \begin{equation}
        \begin{cases}
            Z_{s,h,w} & \text{if } s=1 \\
            Z_{s,h,w} = Z_{s-1,h,w} & \text{otherwise}
        \end{cases}, \forall h\in\mathcal{H}, w\in\mathcal{W}
    \end{equation}

\subsubsection{Hot Water Tank Robust Linking Constraint}
    Hot water tank installation decisions are the same for robust seasons as they are for normal seasons.
    \begin{equation}
            Z_{s_r,h,w} = Z_{s=1,h,w}, \forall s_r\in\mathcal{SR}, h\in\mathcal{H}, w\in\mathcal{W}
    \end{equation}

    \end{appendices}

\end{document}